\documentstyle[11pt,a4wide]{article}
\input amssym.def
\input amssym

\font \tenmsb=msbm10 scaled \magstep 1
\font \sevenmsb=msbm7 scaled \magstep 1
\font \fivemsb=msbm5 scaled \magstep 1
\newfam \msbfam
\textfont \msbfam=\tenmsb
\scriptfont \msbfam=\sevenmsb
\scriptscriptfont \msbfam=\fivemsb
\def \Bbb#1{\fam \msbfam \relax#1}

\font \teneufm=eufm10 scaled \magstep 1
\font \seveneufm=eufm7 scaled \magstep 1
\font \fiveeufm=eufm5 scaled \magstep 1
\newfam \eufmfam
\textfont \eufmfam=\teneufm
\scriptfont \eufmfam=\seveneufm
\scriptscriptfont \eufmfam=\fiveeufm
\def \frak#1{{\fam \eufmfam \relax#1}}

\makeatletter
\@addtoreset{equation}{section}
\makeatother

\title{\bf QUANTUM MATRIX BALL: DIFFERENTIAL AND INTEGRAL CALCULI}

\author{\sl D. Shklyarov \and \sl S. Sinel'shchikov
\and \sl L. Vaksman}

\date{\tt Institute for Low Temperature Physics \& Engineering\\
National Academy of Sciences of Ukraine}

\newtheorem{theorem}{Theorem}[section]
\newtheorem{lemma}[theorem]{Lemma}
\newtheorem{proposition}[theorem]{Proposition}
\newtheorem{corollary}[theorem]{Corollary}

\begin{document}
\large
\maketitle

\bigskip

\section{Introduction}

 The first step in studying q-analogues of irreducible bounded symmetric
domains was made in \cite{SV2}. This work considers the simplest among those
q-analogues, the quantum matrix balls. In this special case, we present here
the proofs of the main results formulated in \cite{SV1, SV2}, and produce,
in particular, an explicit formula for the positive invariant integral, and
thus prove its existence.

 The initial six sections of this work use ${\Bbb C}(q^{1/s})$as the ground
field, the field of rational function of a single indeterminate $q^{1/s}$,
with $s$ being some natural number. The subsequent sections already assume
$q$ to be a number $q \in(0,1)$, and use ${\Bbb C}$ as a ground field.

 We assume a knowledge of the basic notions of quantum group theory
\cite{D}, and, in particular, the notion of a universal R-matrix introduced
by V. Drinfeld. Some of the general properties of a universal R-matrix to be
alluded below, could be easily deduced from the explicit formula for R.
This very well known multiplicative formula for a universal R-matrix is
presented in Appendix 1.

\bigskip

\section{The covariant algebra ${\Bbb C}[{\bf Mat}_{mn}]_q$}

 Recall the definition of the quantum universal enveloping algebra
\index{quantum universal enveloping algebra} $U_q \frak{sl}_N$, introduced
by V. Drinfeld and M. Jimbo. Let $(a_{ij})_{i,j=1,\ldots,N-1}$ be the Cartan
matrix given by
\begin{equation}\label{Cm}
a_{ij}=\left \{\begin{array}{ccl}2 &,& i-j=0 \\
                                 -1 &,& |i-j|=1 \\
                                 0 &,& {\rm otherwise} \end{array}\right..
\end{equation}

The algebra $U_q \frak{sl}_N$ is determined by the generators $E_i$, $F_i$,
$K_i$, $K_i^{-1}$, $i=1,\ldots,N-1$, and the relations

\vbox{
$$K_iK_j=K_jK_i,\quad K_iK_i^{-1}=K_i^{-1}K_i=1,\quad
K_iE_j=q^{a_{ij}}E_jK_i,\quad K_iF_j=q^{-a_{ij}}F_jK_i$$
$$E_iF_j-F_jE_i=\delta_{ij}(K_i-K_I^{-1})/(q-q^{-1})$$
\begin{equation}
E_i^2E_j-(q+q^{-1})E_iE_jE_i+E_jE_i^2=0,\qquad |i-j|=1
\end{equation}
$$F_i^2F_j-(q+q^{-1})F_iF_jF_i+F_jF_i^2=0,\qquad |i-j|=1$$
$$[E_i,E_j]=[F_i,F_j]=0,\qquad |i-j|\ne 1.$$
}

 The comultiplication $\Delta$, the antipode $S$, and the counit
$\varepsilon$ are determined by
\begin{equation}
\Delta(E_i)=E_i \otimes 1+K_i \otimes E_i,\quad \Delta(F_i)=F_i \otimes
K_i^{-1}+1 \otimes F_i,\quad \Delta(K_i)=K_i \otimes K_i,
\end{equation}
\begin{equation}
S(E_i)=-K_i^{-1}E_i,\qquad S(F_i)=-F_iK_i,\qquad S(K_i)=K_i^{-1},
\end{equation}
$$\varepsilon(E_i)=\varepsilon(F_i)=0,\qquad \varepsilon(K_i)=1.$$

 We consider in the sequel only $U_q \frak{sl}_N$-modules of the form
$V=\bigoplus \limits_{\mu \in{\Bbb Z}^{N-1}}V_ \mu$, with
$\mu=(\mu_1,\ldots,\mu_{N-1})$, $V_ \mu=\{v \in
V|\,K_iv=q^{\mu_i}v,\,i=1,\ldots,N-1 \}$. This agreement allows one to
introduce the linear operators $H_j$, $X_j^{\pm}$, $j=1,\ldots,N-1$, by
setting up
\begin{equation}
H_jv=\mu_jv,\quad v \in V_ \mu, \qquad E_j=X_j^+q^{{1 \over2}H_j},\qquad
F_j=q^{-{1 \over 2}H_j}X_j^-.
\end{equation}

 Note that the classical universal enveloping algebra \index{classical
universal enveloping algebra} can be derived from $U_q \frak{sl}_N$ via the
substitution
\begin{equation}q=e^{-h/2},\qquad K_i=e^{-hH_i/2}\end{equation}
and the subsequent passage to a limit as $h \to 0$.

 Turn to a construction of q-analogue of the matrix space ${\rm Mat}_{mn}$,
$m,n \in{\Bbb N}$. Everywhere in the sequel $N=m+n$. We follow \cite{SV2}
in equipping all the $U_q \frak{sl}_N$-modules with the grading ${\rm
deg}\,v=j \:\Leftrightarrow \:  H_0v=2jv$, \index{grading for $U_q
\frak{sl}_N$-modules} where $$H_0={2 \over m+n}\left(m
\sum_{j=1}^{n-1}jH_j+n \sum_{j=1}^{m-1}jH_{N-j}+mnH_n \right).$$ (The
coefficients in the latter identity are chosen so that $[H_0,X_n^\pm]=\pm
2X_n^\pm$, $[H_0,X_j^\pm]=0$ for $j \ne n$.)

 In what follows $V_k$ will stand for the homogeneous components of a graded
vector space $V$, and $V^*$ for the dual graded vector space: $(V^*)_{-k}
\stackrel{\rm def}{=}(V_k)^*$.

 Remind some notions of the theory of Hopf algebras \cite{CP}.

 Let $A$ be a Hopf algebra and $F$ an algebra equipped also by a structure
of $A$-module. $F$ is said to be an $A$-module (covariant) algebra
\index{$A$-module algebra} \index{covariant algebra} if the multiplication
$F \otimes F \to F$, $f_1 \otimes f_2 \mapsto f_1f_2$, is a morphism of
$A$-modules. In the case of a unital algebra $F$, an additional assumption
is introduced that the embedding ${\Bbb C}\hookrightarrow F$, $1 \mapsto 1$,
is a morphism of $A$-modules. A duality argument allows one also to
introduce a notion of $A^{\rm op}$-module (covariant) coalgebra
\index{$A^{\rm op}$-module coalgebra} \index{covariant coalgebra}.

 The notion of a covariant (bi-)module over a covariant algebra and a
covariant (bi-~)comodule over a covariant coalgebra are introduced in a
similar way.

 Let $\lambda=(\lambda_1,\ldots,\lambda_{N-1})\in{\Bbb Z}^{N-1}$, $\lambda_j
\ge 0$ for $j \ne n$. Consider a generalized Verma module
$V_-(\lambda)$ \index{generalized Verma module}. It is a $U_q
\frak{sl}_N$-module with a single generator $v_-(\lambda)$ and the relations
$E_iv_-(\lambda)=0$, $K_i^{\pm 1}v_-(\lambda)=q^{\pm
\lambda_i}v_-(\lambda)$, $i=1,\ldots,N-1$,
$F_j^{\lambda_j+1}v_-(\lambda)=0$, $j \ne n$.

 The $U_q \frak{sl}_N$-module $V_-(0)$ will be equipped with a structure of
covariant coalgebra: $\Delta_-:v_-(0)\mapsto v_-(0) \otimes v_-(0)$, and the
$U_q \frak{sl}_N$-module $V_-(\lambda)$ with a structure of a covariant
bicomodule: $\Delta_-^L:v_-(\lambda)\mapsto v_-(0)\otimes v_-(\lambda)$;
$\Delta_-^R:v_-(\lambda)\mapsto v_-(\lambda)\otimes v_-(0)$.

 In our work \cite{SV2}, a dual algebra ${\Bbb C}[{\rm Mat}_{mn}]_q$ to
covariant coalgebra $V_-(0)$ was considered, together with covariant
bimodules dual to bicomodules $V_-(\lambda)$. (Actually the notation ${\Bbb
C}[\frak{g}_{-1}]_q$ was implicit in \cite{SV2} instead of ${\Bbb C}[{\rm
Mat}_{mn}]_q$, since the exposition of that work was not restricted to the
special case of matrix balls).

 The principal purpose of this section is to describe ${\Bbb C}[{\rm
Mat}_{mn}]_q$ in terms of generators and relations.

 Consider the Hopf subalgebra $U_q \frak{sl}_n \subset U_q \frak{sl}_N$
generated by $E_i$, $F_i$, $K_i^{\pm 1}$, $i=1,2,\ldots,n-1$,
and the Hopf subalgebra $U_q \frak{sl}_m \subset U_q \frak{sl}_N$ generated
by $E_{n+i}$, $F_{n+i}$, $K_{n+i}^{\pm 1}$, $i=1,2,\ldots,m-1$. It follows
from the definitions that the homogeneous component ${\Bbb C}[{\rm
Mat}_{mn}]_{q,1}=\{f \in{\Bbb C}[{\rm Mat}_{mn}]_q|\:{\rm deg}\,f=1 \}$ is a
$U_q \frak{sl}_n \otimes U_q \frak{sl}_m$-module. Prove that this module
splits into the tensor product of a $U_q \frak{sl}_n$-module related to the
vector representation and a $U_q \frak{sl}_m$-module related to the covector
representation.

 Consider the $U_q \frak{sl}_n$-module $U$ and the $U_q \frak{sl}_m$-module
$V$, determined in the bases $\{u_a \}_{a=1,\ldots,n}$, $\{v^\alpha
\}_{\alpha=1,\ldots,m}$ by
$$X_i^+u_a=\left \{\begin{array}{ccl}u_{a-1} &,& a=i+1 \\ 0&,&{\rm
otherwise}\end{array}\right.;\qquad X_{n+i}^+v^\alpha=\left
\{\begin{array}{ccl}v^{\alpha-1} &,& \alpha=m-i+1 \\ 0&,&{\rm
otherwise}\end{array}\right.$$
$$X_i^-u_a=\left \{\begin{array}{ccl}u_{a+1} &,& a=i \\ 0 &,& {\rm
otherwise}\end{array}\right.;\qquad X_{n+i}^-v^\alpha=\left
\{\begin{array}{ccl}v^{\alpha+1} &,& \alpha=m-i \\ 0 &,& {\rm
otherwise}\end{array}\right.$$
$$H_iu_a=\left \{\begin{array}{ccl}u_a &,& a=i \\ -u_a &,& a=i+1 \\ 0&,&{\rm
otherwise}\end{array}\right.;\qquad H_{n+i}v^\alpha=\left
\{\begin{array}{ccl}v^\alpha &,& \alpha=m-i \\ -v^\alpha &,& \alpha=m-i+1 \\
0 &,& {\rm otherwise}\end{array}\right..$$

\medskip

\begin{proposition}\label{zaalph} There exists a unique collection
$\{z_a^\alpha \}_{a=1,\ldots,n;\:\alpha=1,\ldots,m}$, of elements of ${\Bbb
C}[{\rm Mat}_{mn}]_{q,1}$ such that the map $i:u_a \otimes v^\alpha \mapsto
z_a^\alpha$, $a=1,\ldots,n;\:\alpha=1,\ldots,m$ admits an extension up to an
isomorphism of $U_q \frak{sl}_n \otimes U_q \frak{sl}_m$-modules $i:U
\otimes V \mapsto{\Bbb C}[{\rm Mat}_{mn}]_{q,1}$, and $F_nz_n^m=q^{1/2}$.
\end{proposition}

\smallskip

 {\bf Proof.} Consider the maximum length elements for the permutation group
$S_N$ and for its subgroup $S_n \times S_m$
$$w_0=(N,N-1,\ldots,2,1),\qquad w_0'=(n,n-1,\ldots,1,N,N-1,\ldots,n+1).$$
Impose the notation $M=N(N-1)/2$, $M'=M-mn$, $s_j=(j,j+1)$. Let
$w_0=s_{i_1}s_{i_2}\ldots s_{i_M}$ be such a reduced decomposition for $w_0$
that $s_{i_1}s_{i_2}\ldots s_{i_{M'}}=w_0'$. Consider the Hopf subalgebra
$U_q \frak{N}_-\subset U_q \frak{sl}_N$ generated by $\{F_j
\}_{j=1,\ldots,N-1}$, and the base $\left
\{\widetilde{F}_{\beta_M}^{k_M}\widetilde{F}_{\beta_{M-1}}^{k_{M-1}} \ldots
\widetilde{F}_{\beta_1}^{k_1}\right \}_{k_1,\ldots,k_M \in{\Bbb Z}_+}$ in
the vector space $U_q \frak{N}_-$ associated to the above reduced
decomposition.

 The reader is referred to the Appendix 1 for a description of this base,
together with the associated base of the graded vector space $V_-(0)$:
$\left \{\widetilde{F}_{\beta_M}^{k_M}\widetilde{F}_{\beta_{M-1}}^{k_{M-1}}
\ldots \widetilde{F}_{\beta_{M'+1}}^{k_{M'+1}}v_-(0)\right \}$, with
$(k_{M'+1},k_{M'+2},\ldots,k_M)\in{\Bbb Z}_+^{mn}$. Hence, the
dimensionalities of the weight subspaces are just the same as in the
classical $(q=1)$ case. In particular,
\begin{equation}\label{dim}
{\rm dim}\,V_-(0)_{-k}\stackrel{\rm def}{=}{\rm dim}\{v|\:H_0v=-2kv
\}=\left({mn+k-1 \atop k}\right).
\end{equation}
Note that ${\rm dim}\,V_-(0)_{-1}=mn$, and $v'=F_nv_-(0)$ is a non-zero
primitive vector:
$$E_jv'=0,\qquad H_jv'=\left \{\begin{array}{ccl}-2v'&,&j=n \\ v'&,&|j-n|=1
\\ 0&,&|j-n|>1 \end{array}\right.,\qquad j=1,\ldots,N-1.$$
Hence, the $U_q \frak{sl}_n \otimes U_q \frak{sl}_m$-module $U \otimes V$ is
isomorphic to the $U_q \frak{sl}_n \otimes U_q \frak{sl}_m$-module
$(V_-(0)_{-1})^*\simeq{\Bbb C}[{\rm Mat}_{mn}]_{q,1}$. Of course, the
isomorphism $i:U \otimes V \to{\Bbb C}[{\rm Mat}_{mn}]_{q,1}$ is unique up
to a multiple from the ground field, and the elements $z_a^\alpha=i(u_a
\otimes v^\alpha)$, $a=1,\ldots,n$, $\alpha=1,\ldots,m$, satisfy all the
requirements of our proposition, except, possibly, the last property
$F_nz_n^m=q^{1/2}$. One can readily choose the above multiple in the
definition of $i$, which provides this property unless $F_nz_n^m=0$. In the
latter case one has $F_n(E_{i_1}^{k_1}E_{i_2}^{k_2}\ldots
E_{i_l}^{k_l}z_n^m)=0$ for all $i_1,\ldots,i_l$ different from $n$ and all
$k_1, k_2,\ldots,k_l \in{\Bbb Z}_+$. Hence $F_n{\Bbb C}[{\rm
Mat}_{mn}]_{q,1}=0$, and thus $F_nv_-(0)=0$. That is, ${\rm dim}\,V_-(0)=1$.
On the other hand, it follows from (\ref{dim}) that ${\rm
dim}\,V_-(0)=\infty$. This contradiction shows that $F_nz_n^m \ne 0$. \hfill
$\Box$

\medskip

\begin{proposition}\label{gen} $z_a^\alpha$, $a=1,\ldots,n$,
$\alpha=1,\ldots,m$, generate the algebra ${\Bbb C}[{\rm Mat}_{mn}]_q$.
\end{proposition}

\smallskip

 {\bf Proof.} By a virtue of (\ref{dim}), it suffices to prove that for all
$k \in{\Bbb Z}_+$, one can choose $\displaystyle \left({mn+k-1 \atop
k}\right)$ linear independent vectors among the monomials
$z_{a_1}^{\alpha_1}z_{a_2}^{\alpha_2}\ldots z_{a_k}^{\alpha_k}\in{\Bbb
C}[{\rm Mat}_{mn}]_{q,k}$. An application of the standard argument (see
\cite[chapter 5]{Ja}) reduces this statement to its classical analogue.

 Consider the ring $A={\Bbb C}[q^{1/s},q^{-1/s}]$ and the $A$-algebra $U_A$
generated by the elements $E_i$, $F_i$, $K_i^{\pm 1}$,
$L_i={\textstyle K_i-K_i^{-1}\over \textstyle q-q^{-1}}$. This is a Hopf
algebra:  $\Delta(L_i)=L_i \otimes K_i+K_i^{-1}\otimes L_i$, $S(L_i)=-L_i$,
$\varepsilon(L_i)=0$, $i=1,\ldots,N-1$. Let $V_A=U_Av_-(0)$, and $F_A
\subset{\Bbb C}[{\rm Mat}_{mn}]_q$ be the minimal $A$-module which contains
all the monomials $z_{a_1}^{\alpha_1}z_{a_2}^{\alpha_2}\ldots
z_{a_k}^{\alpha_k}$. It follows from the definitions that the value of a
linear functional $z_a^\alpha$ on a vector $v \in V_A$ is in $A$. Hence, a
similar statement is also valid for all the monomials
$z_{a_1}^{\alpha_1}z_{a_2}^{\alpha_2}\ldots z_{a_k}^{\alpha_k}$, $k \in{\Bbb
Z}_+$, and thus for all $f \in F_A$. By means of a specialization $q=1$ we
get (see \cite[chapter 5]{Ja}, \cite{CP}):  $$F_A \to{\Bbb
C}[z_1^1,z_1^2,\ldots,z_n^m],\qquad V_A \to{\Bbb C}\left[{\partial \over
\partial z_1^1},{\partial \over \partial z_1^2},\ldots,{\partial \over
\partial z_n^m}\right].$$ What remains is to apply the non-degeneracy of the
natural pairing for the graded vector spaces ${\Bbb C}[z_1^1,\ldots,z_n^m]$,
${\Bbb C}\left[{\textstyle \partial \over \textstyle \partial
z_1^1},\ldots,{\textstyle \partial \over \textstyle \partial z_n^m}\right]$
and the fact that the dimensionalities of the corresponding homogeneous
components are in both cases $\left({\textstyle mn+k-1 \atop \textstyle
k}\right)$. \hfill $\Box$

\medskip

\begin{proposition}\label{zcr}
\begin{equation}\label{zcr1}
z_{a_1}^{\alpha_1}z_{a_2}^{\alpha_2}-qz_{a_2}^{\alpha_2}z_{a_1}^{\alpha_1}=0,
\qquad a_1=a_2 \quad \&\quad \alpha_1<\alpha_2 \qquad{\rm or}\qquad a_1<a_2
\quad \&\quad \alpha_1=\alpha_2,
\end{equation}
\begin{equation}\label{zcr2}
z_{a_1}^{\alpha_1}z_{a_2}^{\alpha_2}-z_{a_2}^{\alpha_2}z_{a_1}^{\alpha_1}=0,
\qquad \alpha_1<\alpha_2 \quad \&\quad a_1>a_2,
\end{equation}
\begin{equation}\label{zcr3}
z_{a_1}^{\alpha_1}z_{a_2}^{\alpha_2}-z_{a_2}^{\alpha_2}z_{a_1}^{\alpha_1}=
(q-q^{-1})z_{a_1}^{\alpha_2}z_{a_2}^{\alpha_1},\qquad \alpha_1<\alpha_2
\quad \&\quad a_1<a_2.
\end{equation}
\end{proposition}

\smallskip

 {\bf Proof.} The validity of (\ref{zcr1}) -- (\ref{zcr3}) follows from
their validity at the 'classical limit $q=1$' and the $U_q \frak{sl}_n
\otimes U_q \frak{sl}_m$-invariance of the associated subspace in ${\Bbb
C}[{\rm Mat}_{mn}]_{q,1}^{\otimes 2}$. Let $M={\Bbb C}[{\rm
Mat}_{mn}]_{q,1}$ and $M_A \subset M$ be the $A$-module generated by
$\{z_a^\alpha \}$, $a=1,\ldots,n$, $\alpha=1,\ldots,m$. Remind that $M$ is a
module over the Hopf algebra $U_q \frak{sl}_n \otimes U_q \frak{sl}_m$. By a
virtue of proposition \ref{zaalph}, the vector space $M^{\otimes 2}$ admits
a decomposition as a sum of four simple pairwise non-isomorphic $U_q
\frak{sl}_n \otimes U_q \frak{sl}_m$-modules. A similar decomposition is
also valid for $M_A \otimes M_A$, where a specialization at $q=1$ leads to
four pairwise non-isomorphic $U \frak{sl}_n \otimes U
\frak{sl}_m$-modules. By misuse of language, one can say that each submodule
of the $U_q \frak{sl}_n \otimes U_q \frak{sl}_m$-module $M^{\otimes 2}$ is
unambiguously determined by its specialization at $q=1$.  Consider two such
submodules. The first $U_q \frak{sl}_n \otimes U_q \frak{sl}_m$-submodule is
the kernel of the multiplication operator ${\Bbb C}[{\rm
Mat}_{mn}]_{q,1}^{\otimes 2}\to{\Bbb C}[{\rm Mat}_{mn}]_{q,2}$, $f_1 \otimes
f_2 \mapsto f_1f_2$. Another $U_q \frak{sl}_n \otimes U_q
\frak{sl}_m$-submodule is the linear span of the elements given by the left
hand sides of (\ref{zcr1}) -- (\ref{zcr3}). Their specializations at $q=1$
coincide, and hence the $U_q \frak{sl}_n \otimes U_q \frak{sl}_m$-submodules
themselves are the same. \hfill $\Box$

\medskip

\begin{proposition} The relation list (\ref{zcr1}) -- (\ref{zcr3}) is
complete.
\end{proposition}

\smallskip

 {\bf Proof.} Consider a graded unital algebra $F$ determined by degree 1
generators $\{u_a^\alpha \}$, $a=1,\ldots,n$, $\alpha=1,\ldots,m$, and the
relations (\ref{zcr1}) -- (\ref{zcr3}) with the letter '$z$' being replaced
by the latter '$u$'. It is an easy exercise to compute the dimensionalities
of the homogeneous components $F^{(k)}=\{f \in F|\;{\rm deg}\,f=k \}$.
Specifically,
\begin{equation}\label{dim1} {\rm dim}\,F^{(k)}=\left({mn+k-1 \atop
k}\right).
\end{equation}
By a virtue of proposition \ref{zcr}, the map $u_a^\alpha \mapsto
z_a^\alpha$, $a=1,\ldots,n$, $\alpha=1,\ldots,m$, admits an extension up to
a homomorphism $F \to{\Bbb C}[{\rm Mat}_{mn}]_q$. It follows from
proposition \ref{gen} that this homomorphism is onto. What remains is to
apply the relations (\ref{dim}), (\ref{dim1}) to establish the coincidence
of the dimensionalities of the graded components:
$${\rm dim}\,{\Bbb C}[{\rm Mat}_{mn}]_{q,k}={\rm dim}\,F^{(k)},\qquad k
\in{\Bbb Z}_+.\eqno \Box$$

 To conclude, note that the commutation relations (\ref{zcr1}) --
(\ref{zcr3}) were used in a different context by a large number of authors
\cite{CP}.

\bigskip

\section{Differential calculus}

 With \cite{SV2} as a background, we describe a differential calculus on the
quantum space of matrices. An advantage of our approach is that it discovers
an additional surprising symmetry of the standard bicovariant differential
calculus on this quantum space \cite{SV1}.

 Consider the vector $v'=F_nv_-(0)$. It follows from the definitions that
\begin{equation}\label{v'1}
E_jv'=0,\qquad H_jv'=-a_{jn}v',\qquad j=1,\ldots,N-1,
\end{equation}
\begin{equation}\label{v'2} F_i^{-a_{in}+1}v'=0,\qquad i \ne n.
\end{equation}
Here $(a_{ij})$ is a Cartan matrix (\ref{Cm}). Now (\ref{v'1}), (\ref{v'2})
imply

\medskip

\begin{proposition} Consider the $U_q \frak{sl}_N$-module $V_-(\lambda')$
with the highest weight $\lambda'=(-a_{1n},-a_{2n},\ldots,-a_{N-1,n})$. The
map $v_-(\lambda')\mapsto F_nv_-(0)$ admits a unique extension up to a
morphism $\delta_-:V_-(\lambda')\to V_-(0)$ of $U_q \frak{sl}_N$-modules.
\end{proposition}

\medskip

 The graded vector space $\bigwedge^1({\rm Mat}_{mn})_q$ dual to
$V_-(\lambda')$ is a covariant bimodule over ${\Bbb C}[{\rm Mat}_{mn}]_q$.
The adjoint to $\delta_-$ operator $d:{\Bbb C}[{\rm Mat}_{mn}]_q \to
\bigwedge^1({\rm Mat}_{mn})_q$ is called a differential. It follows from the
definitions (see \cite{SV2}) that
$$d(f_1f_2)=df_1 \cdot f_2+f_1 \cdot df_2,\qquad f_1,f_2 \in{\Bbb C}[{\rm
Mat}_{mn}]_q.$$

 Describe the ${\Bbb C}[{\rm Mat}_{mn}]_q$-bimodule $\bigwedge^1({\rm
Mat}_{mn})_q$ in terms of generators and relations. Remind that ${\rm
deg}\,\omega=j \Leftrightarrow H_0 \omega=2j \omega$. Let $\bigwedge^1({\rm
Mat}_{mn})_{q,1}=\{\omega \in \bigwedge^1({\rm Mat}_{mn})_q|\;{\rm
deg}\,\omega=1 \}$.

\medskip

\begin{lemma} $dz_a^\alpha$, $\alpha=1,\ldots,m$, $a=1,\ldots,n$, constitute
a base of the vector space $\bigwedge^1({\rm Mat}_{mn})_{q,1}$.
\end{lemma}

\smallskip

 {\bf Proof.} Since $z_a^\alpha$, $\alpha=1,\ldots,m$, $a=1,\ldots,n$, form
a base of the vector space ${\Bbb C}[{\rm Mat}_{mn}]_{q,1}$, it suffices to
prove that the linear map
$$d:{\Bbb C}[{\rm Mat}_{mn}]_{q,1}\to \bigwedge \nolimits^1({\rm
Mat}_{mn})_{q,1}$$
is one-to-one. Consider the adjoint linear operator
$$\delta_-:V_-(\lambda')_{-1}\to V_-(0)_{-1}.$$
It follows from the definition of the $U_q \frak{sl}_N$-module
$V_-(\lambda')$ that the $U_q \frak{sl}_n \otimes U_q \frak{sl}_m$-module
$V_-(\lambda')_{-1}$
is simple. One can easily deduce from proposition A1.2 that it is
determined by the same relations as the $U_q \frak{sl}_n \otimes U_q
\frak{sl}_m$-module $U^*\otimes V^*$. It was also shown in proposition
\ref{gen} that the $U_q \frak{sl}_n \otimes U_q \frak{sl}_m$-module
$V_-(0)_{-1}$ is simple as well. On the other hand,
$\delta_-|_{V_-(\lambda')_{-1}}$ is a non-zero morphism of $U_q
\frak{sl}_n \otimes U_q \frak{sl}_m$-modules:
$\delta_-v_-(\lambda')=F_nv_-(0)\ne 0$. Hence the restrictions of
$\delta_-$ and $d$ onto the corresponding homogeneous components are
one-to-one. \hfill $\Box$

\medskip

\begin{lemma} The ${\Bbb C}[{\rm Mat}_{mn}]_q$-bimodule $\bigwedge^1({\rm
Mat}_{mn})_q$ is a free left ${\Bbb C}[{\rm Mat}_{mn}]_q$-module:
$$\bigwedge \nolimits^1({\rm Mat}_{mn})_q=\bigoplus_{\alpha=1}^m
\bigoplus_{a=1}^n{\Bbb C}[{\rm Mat}_{mn}]_qdz_a^\alpha,$$ and a free right
${\Bbb C}[{\rm Mat}_{mn}]_q$-module:
$$\bigwedge \nolimits^1({\rm Mat}_{mn})_q=\bigoplus_{\alpha=1}^m
\bigoplus_{a=1}^n dz_a^\alpha{\Bbb C}[{\rm Mat}_{mn}]_q.$$
\end{lemma}

\smallskip

 {\bf Proof.} We are about to prove that the maps
$${\Bbb C}[{\rm Mat}_{mn}]_q \otimes \bigwedge \nolimits^1({\rm
Mat}_{mn})_{q,1}\to \bigwedge \nolimits^1({\rm Mat}_{mn})_q,\qquad f \otimes
\omega \mapsto f \omega,$$
$$\bigwedge \nolimits^1({\rm Mat}_{mn})_{q,1}\otimes{\Bbb C}[{\rm
Mat}_{mn}]_q \to \bigwedge \nolimits^1({\rm Mat}_{mn})_q, \qquad \omega
\otimes f \mapsto \omega f$$
are one-to-one. Their injectivity can be easily derived from a similar
result in the case $q=1$ (cf. the proof of proposition \ref{gen}). What
remains is to use the coincidence of the dimensionalities of homogeneous
components of the graded vector spaces $\bigwedge^1({\rm Mat}_{mn})_q$,
${\Bbb C}[{\rm Mat}_{mn}]_q \otimes \bigwedge^1({\rm Mat}_{mn})_{q,1}$,
$\bigwedge \nolimits^1({\rm Mat}_{mn})_{q,1}\otimes{\Bbb C}[{\rm
Mat}_{mn}]_q$. (The dimensionalities of homogeneous components ${\Bbb
C}[{\rm Mat}_{mn}]_{q,k}$, $k \in{\Bbb Z}_+$, were computed before using a
basis in $V_-(0)$ formed by homogeneous elements. The dimensionalities of
homogeneous components of $\bigwedge^1({\rm Mat}_{mn})_q$ could be found in
a similar way: a basis in $V_-(\lambda')$ could be constructed via an
application of an appropriate reduced decomposition of the complete
permutation $w_0 \in S_N$, together with the associated basis in $U_q
\frak{sl}_N$ (see Appendix 1).) \hfill $\Box$

\medskip

 Of course, the covariant ${\Bbb C}[{\rm Mat}_{mn}]_q$-{\sl bimodule}
$\bigwedge^1({\rm Mat}_{mn})_q$ is not free. The elements $dz_a^\alpha$ are
its generators. Find a complete list of relations.

 Let $U_q \frak{sl}_N^{\rm op}$ be the Hopf algebra which differs from $U_q
\frak{sl}_N$ by a replacement of its comultiplication $\Delta$ with an
opposite one $\Delta^{\rm op}$. The structure of a ${\Bbb C}[{\rm
Mat}_{mn}]_q$-bimodule $\bigwedge^1({\rm Mat}_{mn})_q$ has been defined in
\cite{SV2} via an application of a duality argument and the following
morphisms in the category of modules over the Hopf algebra $U_q
\frak{sl}_N^{\rm op}$:
$$\Delta_-^L:V_-(\lambda')\to V_-(0)\otimes V_-(\lambda');\qquad
\Delta_-^L:v_-(\lambda')\mapsto v_-(0)\otimes v_-(\lambda');$$
$$\Delta_-^R:V_-(\lambda')\to V_-(\lambda')\otimes V_-(0);\qquad
\Delta_-^R:v_-(\lambda')\mapsto v_-(\lambda')\otimes v_-(0).$$

 Let $P:V_-(\lambda')\otimes V_-(0)\to V_-(0)\otimes V_-(\lambda')$ be the
ordinary flip of tensor multiples: $P(v'\otimes v'')=v''\otimes v'$.  Define
an operator
$\widetilde{R}_{V_-(\lambda')\,V_-(0)}:V_-(\lambda')\otimes V_-(0)\to
V_-(0)\otimes V_-(\lambda')$ via the universal R-matrix (see Appendix 1) by
$\widetilde{R}_{V_-(\lambda')\,V_-(0)}=R_{V_-(0)\,V_-(\lambda')}P$.

\medskip

\begin{lemma}\label{drdl}
$\widetilde{R}_{V_-(\lambda')\,V_-(0)}\Delta_-^R=\Delta_-^L$.
\end{lemma}

\smallskip

 {\bf Proof.} It is well known that the operator $P \cdot
R_{V_-(\lambda')\,V_-(0)}$ is a morphism in the category of $U_q
\frak{sl}_N$-modules. Hence, the operator $R_{V_-(0)\,V_-(\lambda')}P$ is a
morphism in the category of $U_q \frak{sl}_N^{\rm op}$-modules. What remains
is to apply the identity
$\widetilde{R}_{V_-(\lambda')\,V_-(0)}v_-(\lambda')\otimes
v_-(0)=v_-(0)\otimes v_-(\lambda')$, which follows from the property (A1.6)
of the universal R-matrix. \hfill $\Box$

\medskip

 The universal R-matrix satisfies the identity $S \otimes S(R)=R$, with $S$
being the antipode of the Hopf algebra in question (see \cite{CP}). Hence,
the adjoint to $\widetilde{R}_{V_-(\lambda')\,V_-(0)}$ operator is of the
form
\begin{equation}
\check{R}_{{\Bbb C}[{\rm Mat}_{mn}]_q \,\wedge^1({\rm Mat}_{mn})_q}=PR_{{\Bbb
C}[{\rm Mat}_{mn}]_q \,\wedge^1({\rm Mat}_{mn})_q}
\end{equation}
with $P:{\Bbb C}[{\rm Mat}_{mn}]_q \otimes \wedge^1({\rm Mat}_{mn})_q \to
\wedge^1({\rm Mat}_{mn})_q \otimes {\Bbb C}[{\rm Mat}_{mn}]_q$ being the
ordinary flip of tensor multiples.

\medskip

\begin{corollary} For all $\alpha,\beta=1,\ldots,m$, $a,b=1,\ldots,n$,
\begin{equation}\label{zdz} z_b^\beta dz_a^\alpha=m_RPR_{{\Bbb C}[{\rm
Mat}_{mn}]_q \,\wedge^1({\rm Mat}_{mn})_q}(z_b^\beta \otimes dz_a^\alpha),
\end{equation}
with $m_R:\wedge^1({\rm Mat}_{mn})_q \otimes{\Bbb C}[{\rm Mat}_{mn}]_q \to
\wedge^1({\rm Mat}_{mn})_q$, $m_R:\omega \otimes f \mapsto \omega f$.
\end{corollary}

\smallskip

 {\bf Proof.} It suffices to pass in the statement of lemma \ref{drdl} to
dual graded vector spaces and to adjoint operators. \hfill $\Box$

\medskip

 Simplify (\ref{zdz}) by computing $R_{{\Bbb C}[{\rm Mat}_{mn}]_q
\,\wedge^1({\rm Mat}_{mn})_q}(z_b^\beta \otimes dz_a^\alpha)$ via an
application of the multiplicative formula for the universal R-matrix.

\medskip

\begin{lemma}\label{orth} $H_0$ is orthogonal to all the vectors $H_j$, $j
\ne n$, with respect to the bilinear invariant scalar product
$(H_i,H_j)=a_{ij}$, $i,j=1,\ldots,N-1$.
\end{lemma}

\smallskip

 {\bf Proof.} The invariant scalar product $(H_0,H_j)$ is given by ${\rm
tr}\,\pi_1(H_0)\pi_1(H_j)$, with $\pi_1$ being the vector representation of
the Lie algebra $\frak{sl}_N$. What remains is to compute this trace using
the standard basis $\{e_j \}_{j=1}^N$ and the relation
$$\pi_1(H_0)e_j={2 \over m+n}\left \{\begin{array}{r|c}me_j & j \le n \\
-ne_j & j>n \end{array}\right..\eqno \Box$$

\medskip

 Consider the $U_q \frak{sl}_n \otimes U_q \frak{sl}_m$-modules $L'={\Bbb
C}[{\rm Mat}_{mn}]_{q,1}$, $L''=\wedge^1({\rm Mat}_{mn})_{q,1}$ (the
homogeneous components of the graded vector spaces ${\Bbb C}[{\rm
Mat}_{mn}]_q$, $\wedge^1({\rm Mat}_{mn})_q$). Let $R_{L'L''}$ stand for the
linear operator in $L'\otimes L''$ determined by the action of the universal
R-matrix of the Hopf algebra $U_q \frak{sl}_n \otimes U_q
\frak{sl}_m$.\footnote{This universal R-matrix is a tensor product of the
universal R-matrices (A1.6) for $U_q \frak{sl}_n$ and $U_q \frak{sl}_m$.}

\medskip

\begin{lemma}\label{rll0} For all $\alpha,\beta=1,\ldots,m$,
$a,b=1,\ldots,n$,
\begin{equation}\label{rll}
R_{{\Bbb C}[{\rm Mat}_{mn}]_q \,\wedge^1({\rm Mat}_{mn})_q}(z_b^\beta
\otimes dz_a^\alpha)={\tt const}\cdot R_{L'L''}(z_b^\beta \otimes
dz_a^\alpha),
\end{equation}
with {\tt const} being independent of $a$, $b$, $\alpha$, $\beta$.
\end{lemma}

\smallskip

 {\bf Proof.} Apply to both sides of (\ref{rll}) the multiplicative formula
for the universal R-matrix (A1.6).  The
'redundant' exponential multiples in the left hand side of the resulting
identity can be omitted since
$$\exp_{q^2}((q^{-1}-q)E_{\beta_j}\otimes F_{\beta_j}
)z_b^\beta \otimes
dz_a^\alpha=z_b^\beta \otimes dz_a^\alpha$$
for all $\alpha,\beta=1,\ldots,m$, $a,b=1,\ldots,n$, $j>{\textstyle
m(m-1)\over \textstyle 2}+{\textstyle n(n-1)\over \textstyle 2}$. What
remains is to compare the multiple $q^{-t_0}$ related to the Hopf algebra
$U_q \frak{sl}_N$ to a similar multiple related to the Hopf subalgebra $U_q
\frak{sl}_n \otimes U_q \frak{sl}_m$. It follows from lemma \ref{orth} and
the description of $t_0$ in terms of the orthogonal basis of the Cartan
subalgebra (see Appendix 1) that their actions on the subspace ${\Bbb C}[{\rm
Mat}_{mn}]_{q,1}\otimes \wedge^1({\rm Mat}_{mn})_{q,1}$ differ only by a
constant multiple. \hfill $\Box$

\medskip

 Let
$$\left.\widehat{R}_{UU}\right.^{b'a'}_{ba}=\left \{\begin{array}{ccl}q^{-1}
&,& a=b=a'=b' \\ 1 &,& a \ne b \quad \&\quad a=a'\quad \&\quad b=b'\\
q^{-1}-q &,& a<b \quad \&\quad a=b' \quad \&\quad b=a'\\ 0 &,& {\rm
otherwise}\end{array}\right.,$$
$$\left.\widehat{R}_{VV}\right.^{\beta'\alpha'}_{\beta \alpha}=\left
\{\begin{array}{ccl}q^{-1} &,& \alpha=\beta=\alpha'=\beta' \\ 1 &,& \alpha
\ne \beta \quad \&\quad \alpha=\alpha'\quad \&\quad \beta=\beta'\\ q^{-1}-q
&,& \alpha<\beta \quad \&\quad \alpha=\beta' \quad \&\quad \beta=\alpha'\\ 0
&,& {\rm otherwise}\end{array}\right..$$

\medskip

\begin{proposition}\label{zdz1} For all $\alpha,\beta=1,\ldots,m$,
$a,b=1,\ldots,n$,
$$z_b^\beta dz_a^\alpha=\sum_{\alpha',\beta'=1}^m \sum_{a',b'=1}^n
\left.\widehat{R}_{VV}\right._{\beta'\alpha'}^{\beta
\alpha}\left.\widehat{R}_{UU}\right.^{b'a'}_{ba}dz_{a'}^{\alpha'}\cdot
z_{b'}^{\beta'}.$$
\end{proposition}

\smallskip

 {\bf Proof.} Consider the operators in $U \otimes U$ and $V \otimes V$
determined by the actions of the universal R-matrices for Hopf algebras $U_q
\frak{sl}_n$ and $U_q \frak{sl}_m$ respectively. It is well known  (see
\cite{D, CP}) that these operators coincide up to constant multiples with
the operators $\widehat{R}_{UU}$, $\widehat{R}_{VV}$ given by the matrices
$\left.\widehat{R}_{UU}\right.^{b'a'}_{ba}$,
$\left.\widehat{R}_{VV}\right._{\beta'\alpha'}^{\beta \alpha }$. Hence,
by virtue of (\ref{zdz}), (\ref{rll}), and proposition \ref{zaalph}, one has
$$z_b^\beta dz_a^\alpha={\tt const}_1 \sum_{\alpha',\beta'=1}^m
\sum_{a',b'=1}^n \left.\widehat{R}_{VV}\right._{\beta'\alpha'}^{\beta
\alpha}\left.\widehat{R}_{UU}\right.^{b'a'}_{ba}dz_{a'}^{\alpha'}\cdot
z_{b'}^{\beta'}.$$
What remains is to prove that ${\tt const}_1=1$. This is due to
$$\langle z_n^mdz_n^m,F_nv_-(\lambda')\rangle=q^{-2}\langle dz_n^m \cdot
z_n^m,F_nv_-(\lambda')\rangle \ne 0.$$
The latter relation could be easily deduced from the definitions (just as it
was done in the special case $m=n=1$ described in details in \cite{SV2}).
\hfill $\Box$

\medskip

 We have described an \underline{order one} differential calculus
\index{order one differential calculus} on the quantum matrix space in terms
of generators and relations. Consider the associated universal
\underline{full} differential calculus \index{full differential calculus}
(see, for instance \cite{SV2}).  Proposition \ref{zdz1} implies

\medskip

\begin{corollary}
$$dz_b^\beta dz_a^\alpha=-\sum_{\alpha',\beta'=1}^m \sum_{a',b'=1}^n
\left.\widehat{R}_{VV}\right._{\beta'\alpha'}^{
\beta \alpha}\left.\widehat{R}_{UU}\right.^{b'a'}_{ba}dz_{a'}^{\alpha'}\cdot
dz_{b'}^{\beta'}.$$
\end{corollary}

\medskip

 The differential algebra $\wedge({\rm Mat}_{mn})_q$ \index{differential
algebra $\wedge({\rm Mat}_{mn})_q$} described here in terms of generators
and relations is well known \cite{CP}. Our approach to its construction made
it possible to discover a hidden symmetry of this differential algebra (see
\cite{SV1}). While producing the covariant algebras ${\Bbb C}[{\rm
Mat}_{mn}]_q$, $\wedge({\rm Mat}_{mn})_q$, the generalized Verma modules
$V_-(\lambda)\ni v_-(\lambda)$ with highest weights were implemented. It was
demonstrated in \cite{SV2} that, after replacing them by the generalized
Verma modules $V_+(\lambda)\ni v_+(\lambda)$ with lowest weights, it is
possible to produce a covariant algebra of 'antiholomorphic polynomials' and
the associated differential algebra $\wedge(\overline{\rm Mat}_{mn})_q$.

\bigskip

\section{Covariant $\bf *$-algebra ${\bf Pol}({\bf Mat}_{mn})_q$}

 Remind \cite{CP} that in the case of involutive algebras the definition of
an $A$-module algebra includes the following compatibility axiom for
involutions:
\begin{equation}\label{*H} (af)^*=(S(a))^*f^*,\qquad a \in A,\;f \in F
\end{equation}

 Let $U_q \frak{su}_{nm}$ \index{$U_q \frak{su}_{nm}$} stand for the
$*$-Hopf algebra $(U_q \frak{sl}_N,*)$ given by
$$(K_j^{\pm 1})^*=K_j^{\pm 1},\qquad E_j^*=\left \{\begin{array}{rcc}K_jF_j
&,& j \ne n \\ -K_jF_j &,& j=n \end{array}\right.,\qquad F_j^*=\left
\{\begin{array}{rcc}E_jK_j^{-1} &,& j \ne n \\ -E_jK_j^{-1} &,&
j=n \end{array}\right.,$$
with $j=1,\ldots,N-1$. In terms of the 'generators' $H_j$, $X_j^{\pm 1}$ (i.
e. for operators from the class of $*$-representations of $U_q
\frak{su}_{nm}$ described in section 2) one has
$$H_j^*=H_j,\qquad (X_j^\pm)^*=\left \{\begin{array}{rcc}X_j^\mp &,& j \ne
n \\ -X_j^\mp &,& j=n \end{array}\right.,\quad j=1,\ldots,N-1.$$

 A standard method of quantum group theory was used in \cite{SV2} to equip
each of the spaces
$${\rm Pol}({\rm Mat}_{mn})_q \stackrel{\rm def}{=}{\Bbb C}[{\rm Mat}_{mn}]_q
\otimes{\Bbb C}[\overline{\rm Mat}_{mn}]_q,\qquad \Omega({\rm Mat}_{mn})_q
\stackrel{\rm def}{=}\bigwedge({\rm Mat}_{mn})_q \otimes
\bigwedge(\overline{\rm Mat}_{mn})_q$$
with a structure of $U_q \frak{sl}_{nm}$-module algebra (covariant algebra).
The subalgebras
$${\Bbb C}[{\rm Mat}_{mn}]_q \otimes 1 \subset{\rm Pol}({\rm
Mat}_{mn})_q,\qquad 1 \otimes{\Bbb C}[\overline{\rm Mat}_{mn}]_q \subset{\rm
Pol}({\rm Mat}_{mn})_q$$
are conjugate ($*:{\Bbb C}[{\rm Mat}_{mn}]_q \to{\Bbb C}[\overline{\rm
Mat}_{mn}]_q$); they are q-analogues of subalgebras of holomorphic and
antiholomorphic polynomials respectively.

 It follows from the definitions of \cite{SV2} and proposition \ref{gen}
that $\{z_a^\alpha \}$, $\alpha=1,\ldots,m$, $a=1,\ldots,n$, generate the
$*$-algebra ${\rm Pol}({\rm Mat}_{mn})_q$, and the complete relation list
consists of (\ref{zcr1}) -- (\ref{zcr3}), together with the following
R-matrix commutation relation (cf. (\ref{zdz})):
\begin{equation}\label{zsz} (z_b^\beta)^*z_a^\alpha=mPR_{{\Bbb
C}[\overline{\rm Mat}_{mn}]_q \,{\Bbb C}[{\rm
Mat}_{mn}]_q}(z_b^\beta)^*\otimes z_a^\alpha,
\end{equation}
with $m:{\rm Pol}({\rm Mat}_{mn})_q^{\otimes 2}\to{\rm Pol}({\rm
Mat}_{mn})_q$, $m:f_1 \otimes f_2 \mapsto f_1f_2$ being the multiplication
in ${\rm Pol}({\rm Mat}_{mn})_q$, $P$ the flip of tensor multiples, and
$R_{{\Bbb C}[\overline{\rm Mat}_{mn}]_q \,{\Bbb C}[{\rm Mat}_{mn}]_q}$ the
linear operator in ${\Bbb C}[\overline{\rm Mat}_{mn}]_q \otimes{\Bbb C}[{\rm
Mat}_{mn}]_q$ determined by the universal R-matrix.

 Simplify the expression $R_{{\Bbb C}[\overline{\rm Mat}_{mn}]_q \,{\Bbb
C}[{\rm Mat}_{mn}]_q}(z_b^\beta)^*\otimes z_a^\alpha$ and thus the right
hand side of (\ref{zsz}).

 Denote by $U_q \frak{su}_n \otimes U_q \frak{su}_m$ \index{$U_q \frak{su}_n
\otimes U_q \frak{su}_m$} the subalgebra of the $*$-Hopf algebra $U_q
\frak{su}_{nm}$ generated by $E_j$, $F_j$, $K_j$, $K_j^{-1}$ with $j \ne n$.

 Now an application of proposition \ref{zaalph} makes it easy to prove the
following

\medskip

\begin{lemma}\label{spi} The sesquilinear form in ${\Bbb
C}[{\rm Mat}_{mn}]_{q,1}$ given by
$(z_a^\alpha,z_b^\beta)=\delta_{ab}\delta^{\alpha \beta}$, $a,b=1,\ldots,n$,
$\alpha,\beta=1,\ldots,m$, is $U_q \frak{su}_n \otimes U_q
\frak{su}_m$-invariant: $(\xi
z_a^\alpha,z_b^\beta)=(z_a^\alpha,\xi^*z_b^\beta)$ for all $\xi \in U_q
\frak{su}_n \otimes U_q \frak{su}_m$, $a,b=1,\ldots,n$,
$\alpha,\beta=1,\ldots,m$.
\end{lemma}

\medskip

 Note that $(z_b^\beta)^*$, $b=1,\ldots,n$, $\beta=1,\ldots,m$, form a base
for the homogeneous component ${\Bbb C}[\overline{\rm Mat}_{mn}]_{q,-1}$ of
the graded vector space ${\Bbb C}[\overline{\rm Mat}_{mn}]_q$.

\medskip

\begin{corollary}\label{mu} The linear functional $\mu$ on ${\Bbb
C}[\overline{\rm Mat}_{mn}]_{q,-1}\otimes{\Bbb C}[{\rm Mat}_{mn}]_{q,1}$
given by $\mu((z_b^\beta)^*\otimes z_a^\alpha)=\delta_{ab} \delta^{\alpha
\beta}$, is invariant (i. e. $\mu(\xi((z_b^\beta)^*\otimes
z_a^\alpha))=\varepsilon(\xi)\mu((z_b^\beta)^*\otimes z_a^\alpha)$ for all
$\xi \in U_q \frak{su}_n \otimes U_q \frak{su}_m$, $a,b=1,\ldots,n$,
$\alpha,\beta=1,\ldots,m$).
\end{corollary}

\smallskip

 {\bf Proof.} Let $L={\Bbb C}[{\rm Mat}_{mn}]_{q,1}$. Consider the
\underline{antimodule} \index{antimodule} $\overline{L}$ which is still $L$
as an Abelian group, but the actions of the ground field and $U_q
\frak{su}_n \otimes U_q \frak{su}_m$ are given by $(\lambda,v)\mapsto
\overline{\lambda}v$, $(\xi,v)\mapsto S(\xi)^*v$, $\xi \in U_q \frak{su}_n
\otimes U_q \frak{su}_m$, $v \in L$. It follows from lemma \ref{spi} that
the linear functional $\overline{L}\otimes L \to {\Bbb C}(q^{1/s})$
corresponding to the sesquilinear form in $L$, is invariant.
The relationship of an invariant integral and an invariant form is discussed,
 for example, in \cite{SSV2}. \hfill $\Box$

\medskip

 Let $L'={\Bbb C}[\overline{\rm Mat}_{mn}]_{q,-1}$, $L''={\Bbb C}[{\rm
Mat}_{mn}]_{q,1}$, and $R_{L'L''}$ is the linear operator in $L'\otimes L''$
given by the action of the universal R-matrix of the Hopf algebra $U_q
\frak{sl}_n \otimes U_q \frak{sl}_m \subset U_q \frak{sl}_N$.

\medskip

\begin{lemma}\label{c1c2} For all $a,b=1,\ldots,n$,
$\alpha,\beta=1,\ldots,m$,
$$R_{{\Bbb C}[\overline{\rm Mat}_{mn}]_q \,{\Bbb C}[{\rm
Mat}_{mn}]_q}((z_b^\beta)^*\otimes z_a^\alpha)={\tt const}_1
\cdot R_{L'L''}((z_b^\beta)^*\otimes z_a^\alpha)+{\tt const}_2 \cdot
\delta_{ab}\delta^{\alpha \beta},$$
with ${\tt const}_1$ and ${\tt const}_2$ being independent of $a$, $b$,
$\alpha$, $\beta$.
\end{lemma}

\smallskip

 {\bf Proof.} Reproduce essentially the proof of lemma \ref{rll0} to
establish the existence of such element ${\tt const}_1$ of the ground field
that for all $a$, $b$, $\alpha$, $\beta$ one has
$$R_{{\Bbb C}[\overline{\rm Mat}_{mn}]_q \,{\Bbb C}[{\rm
Mat}_{mn}]_q}((z_b^\beta)^*\otimes z_a^\alpha)-{\tt const}_1
\cdot R_{L'L''}((z_b^\beta)^*\otimes z_a^\alpha)\in{\Bbb C}[\overline{\rm
Mat}_{mn}]_{q,0} \otimes{\Bbb C}[{\rm Mat}_{mn}]_{q,0}.$$
Thus we get a linear functional on ${\Bbb C}[\overline{\rm Mat}_{mn}]_{q,-1}
\otimes{\Bbb C}[{\rm Mat}_{mn}]_{q,1}$ since
$${\rm dim}({\Bbb C}[\overline{\rm Mat}_{mn}]_{q,0})={\rm dim}({\Bbb C}[{\rm
Mat}_{mn}]_{q,0})=1.$$
In virtue of general properties of a universal R-matrix \cite{CP} this
linear functional belongs to the subspace of $U_q \frak{sl}_n \otimes U_q
\frak{sl}_m$-invariant linear functionals. This subspace is one-dimensional
due to the simplicity of the $U_q \frak{sl}_n \otimes U_q
\frak{sl}_m$-module ${\Bbb C}[{\rm Mat}_{mn}]_{q,1}$. What remains is to
apply corollary \ref{mu}. \hfill $\Box$

\medskip

 We need an explicit form of the operator $R_{L'L''}$.

 Let $*:U \to \overline{U}$, $*:V \to \overline{V}$, be the identical maps
between the above $U_q \frak{sl}_n$-module $U$ and $U_q \frak{sl}_m$-module
$V$ onto the associated antimodules. Let $R_{\overline{U}U}$,
$R_{\overline{V}V}$ stand for the operators in
$\overline{U}\otimes U$, $\overline{V}\otimes V$ respectively, given by the
actions of the universal R-matrices of the Hopf algebras $U_q \frak{sl}_n$
and $U_q \frak{sl}_m$.

\medskip

\begin{lemma}\label{c'c''} For all $a,b=1,\ldots,n$,
$\alpha,\beta=1,\ldots,m$,
$$R_{\overline{U}U}u_b^*\otimes u_a={\tt const}'\cdot \left
\{\begin{array}{ccl}q^{-1}u_b^*\otimes u_a &,& a \ne b \\ u_a^*\otimes
u_a-(q^{-2}-1)\sum \limits_{k>a}u_k^* \otimes u_k &,&a=b
\end{array}\right.,$$
$$R_{\overline{V}V}(v^\beta)^*\otimes v^\alpha={\tt const}''\cdot
\left \{\begin{array}{ccl}q^{-1}(v^\beta)^*\otimes v^\alpha &,& \alpha
\ne \beta \\ (v^\alpha)^*\otimes v^\alpha-(q^{-2}-1)\sum
\limits_{k>\alpha}(v^k)^*\otimes v^k &,& \alpha=\beta \end{array}\right..$$
with ${\tt const}'$, ${\tt const}''$ being independent of $a$, $b$,
$\alpha$, $\beta$.
\end{lemma}

\smallskip

 {\bf Proof.} It suffices to prove the first identity. Consider the linear
operator $PR_{\overline{U}U}:\overline{U}\otimes U \to{U}\otimes
\overline{U}$, with $P$ being the flip of tensor multiples. It follows from
the general properties of the universal R-matrix that this operator is a
morphism of $U_q \frak{sl}_n$-modules. Besides, it follows from (A1.6) that
$PR_{\overline{U}U}u_n^* \otimes u_n={\tt const}'\cdot u_n \otimes u_n^*$
since $u_n$ is the lowest weight vector of the $U_q \frak{sl}_n$-module $U$.

 On the other hand, it is well known (see, for example, \cite{SoV1}) that
the operators defined by the right hand sides of the identities in the
statement of our lemma possess the same properties. What remains is to use
the fact that each morphism of $U_q \frak{sl}_n$-modules
$\overline{U}\otimes U \to{U}\otimes \overline{U}$ which annihilates $u_n^*
\otimes u_n$, is identically zero (this vector does not belong to any of
the two simple components of the $U_q \frak{sl}_n$-module
$\overline{U}\otimes U $, and hence it generates
this module). \hfill $\Box$

\medskip

 Lemmas \ref{c1c2}, \ref{c'c''} allow one to deduce all the relations
between $(z_b^\beta)^*$, $z_a^\alpha$ up to two constants. These will be
computed by means of the following

\medskip

\begin{lemma} $R_{{\Bbb C}[\overline{\rm Mat}_{mn}]_q \,{\Bbb C}[{\rm
Mat}_{mn}]_q}(z_n^m)^*\otimes z_n^m=q^2(z_n^m)^*\otimes z_n^m+1-q^2$.
\end{lemma}

\smallskip

 {\bf Proof.} We are about to apply the explicit formula (A1.6) for the
universal R-matrix.

 Prove that $H_jz_n^m=\left \{\begin{array}{ccl} 2z_n^m &,& j=n \\ -z_n^m
&,& |j-n|=1 \\ 0 &,& {\rm otherwise}\end{array}\right.$. The two latter
relations follow from the definitions of $z_a^\alpha$, see section 2. The
first relation follows from $H_0z_n^m=2z_n^m$:
$$2z_n^m={2 \over m+n}(-m(n-1)-n(m-1))z_n^m+{2mn \over m+n}H_nz_n^m.$$
Hence $z_n^m$, $(z_n^m)^*$ are weight vectors whose weights are $\alpha_n$,
$-\alpha_n$ respectively. Thus, we have
$$t_0((z_n^m)^*\otimes z_n^m)=(-\alpha_n,\alpha_n)(z_n^m)^*\otimes
z_n^m=-2(z_n^m)^*\otimes z_n^m.$$
What remains is to take into account that only q-exponent survives in
(A1.6), and to use the detailed calculations for the case $m=n=1$ given in
\cite{SV2}. \hfill $\Box$

\medskip

\begin{corollary}\label{znmcr} $(z_n^m)^*z_n^m=q^2z_n^m(z_n^m)^*+1-q^2$.
\end{corollary}

\medskip

 Let
$$\left.\widehat{R}_{\overline{U}U}\right._{ba}^{b'a'}=\left
\{\begin{array}{ccl}q^{-1} &,& a \ne b \quad \&\quad b=b'\quad \&\quad
a=a'\\ 1 &,& a=b=a'=b'\\ -(q^{-2}-1) &,& a=b \quad \&\quad a'=b'\quad
\&\quad a'>a \\ 0 &,&{\rm otherwise}\end{array}\right.,$$
$$\left.\widehat{R}_{\overline{V}V}\right._{\beta'\alpha'}^{\beta
\alpha}=\left \{\begin{array}{ccl}q^{-1} &,& \alpha \ne \beta \quad \&\quad
\beta=\beta'\quad \&\quad \alpha=\alpha'\\ 1 &,&
\alpha=\beta=\alpha'=\beta'\\ -(q^{-2}-1) &,& \alpha=\beta \quad \&\quad
\alpha'=\beta'\quad \&\quad \alpha'>\alpha \\ 0 &,&{\rm
otherwise}\end{array}\right..$$

\medskip

\begin{proposition} For all $a,b=1,\ldots,n$, $\alpha,\beta=1,\ldots,m$,
\begin{equation} (z_b^\beta)^*\cdot z_a^\alpha=q^2 \cdot \sum_{a',b'=1}^n
\sum_{\alpha',\beta'=1}^m
\left.\widehat{R}_{\overline{U}U}\right._{ba}^{b'a'}\cdot
\left.\widehat{R}_{\overline{V}V}\right._{\beta'\alpha'}^{\beta
\alpha}z_{a'}^{\alpha'}(z_{b'}^{\beta'})^*+(1-q^2)\delta_{ab}\delta^{\alpha
\beta}.
\end{equation}
\end{proposition}

\smallskip

 {\bf Proof.} The desired commutation relation with indefinite coefficients
instead of $q^2$ and $1-q^2$ follows from lemmas \ref{c1c2}, \ref{c'c''}.
The values of those coefficients can be found via an application of
corollary \ref{znmcr}. \hfill $\Box$

\medskip

 An application of the operators $\partial$, $\overline{\partial}$ (see
\cite{SV2}) yields

\medskip

\begin{corollary} For all $a,b=1,\ldots,n$, $\alpha,\beta=1,\ldots,m$,
$$d(z_b^\beta)^*\cdot z_a^\alpha=q^2 \cdot \sum_{a',b'=1}^n
\sum_{\alpha',\beta'=1}^m
\left.\widehat{R}_{\overline{U}U}\right._{ba}^{b'a'}\cdot
\left.\widehat{R}_{\overline{V}V}\right._{\beta'\alpha'}^{\beta
\alpha}z_{a'}^{\alpha'}d(z_{b'}^{\beta'})^*,$$
$$(z_b^\beta)^*\cdot dz_a^\alpha=q^2 \cdot \sum_{a',b'=1}^n
\sum_{\alpha',\beta'=1}^m
\left.\widehat{R}_{\overline{U}U}\right._{ba}^{b'a'}\cdot
\left.\widehat{R}_{\overline{V}V}\right._{\beta'\alpha'}^{\beta
\alpha}dz_{a'}^{\alpha'}(z_{b'}^{\beta'})^*,$$
$$d(z_b^\beta)^*\cdot dz_a^\alpha=-q^2 \cdot \sum_{a',b'=1}^n
\sum_{\alpha',\beta'=1}^m
\left.\widehat{R}_{\overline{U}U}\right._{ba}^{b'a'}\cdot
\left.\widehat{R}_{\overline{V}V}\right._{\beta'\alpha'}^{\beta
\alpha}dz_{a'}^{\alpha'}d(z_{b'}^{\beta'})^*.$$
\end{corollary}

\bigskip

\section{The quantum group $\bf SL_N$}

 Remind that currently we use ${\Bbb C}(q^{1/s})$, $s \in{\Bbb N}$, as a
ground field. Later on, we shall, keeping the notation, pass to ${\Bbb C}$
as a ground field.

 Consider the Hopf algebra ${\Bbb C}[SL_N]_q$ of regular functions on the
quantum group $SL_N$ (see \cite{D, RTF}). This algebra is determined by its
generators $\{t_{ij}\}_{i,j=1,\ldots,N}$, the commutation relations
analogous to (\ref{zcr1}) -- (\ref{zcr3}), and the relation ${\rm
det}_qT=1$. (Here ${\rm det}_qT$ is a q-determinant \index{q-determinant} of
the matrix $T=(t_{ij})_{i,j=1,\ldots,N}$:
$${\rm det}_qT=\sum_{s \in S_N}(-q)^{l(s)}t_{1s(1)}t_{2s(2)}\ldots
t_{Ns(N)},$$
with $l(s)={\rm card}\{(i,j)|\;i<j \quad \&\quad s(i)>s(j) \}$).
Comultiplication $\Delta$, counit $\varepsilon$, and antipode $S$ are
defined as follows:
$$\Delta(t_{ij})=\sum_kt_{ik}\otimes t_{kj},\qquad
\varepsilon(t_{ij})=\delta_{ij},\qquad S(t_{ij})=(-q)^{i-j}{\rm
det}_qT_{ji}.$$
Here $i,j=1,\ldots,N$, and the matrix $T_{ji}$ is derived from $T$ by
obliterating its $j$-th line and $i$-th column.

 Just as in section 2, we consider the vector representation $\pi_1$ of the
Hopf algebra $U_q \frak{sl}_N$ and the basis $\{u_k \}_{k=1,\ldots,N}$ in
the space of this representation. We also need a well known \cite{CP}
non-degenerate pairing of {\sl Hopf algebras} ${\Bbb C}[SL_N]_q \times U_q
\frak{sl}_N \to{\Bbb C}(q^{1/s})$ in which a pair $(t_{ij},\xi)$,
$i,j=1,\ldots,N$, is sent to the corresponding matrix element of the
operator $\pi_1(\xi)$.

 This pairing is used to equip ${\Bbb C}[SL_N]_q$ with a structure of $U_q
\frak{sl}_N^{\rm op}\otimes U_q \frak{sl}_N$-module algebra as follows:
\begin{equation} \langle(\eta \otimes \xi)f,\zeta \rangle=\langle
f,S(\eta)\zeta \xi \rangle
\end{equation}
for all $f \in{\Bbb C}[SL_N]_q$, $\xi,\eta,\zeta \in U_q \frak{sl}_N$ (see
the definition of $U_q \frak{sl}_N^{\rm op}$ in section 3).

 We have described a q-analogue for the action of $SL_N \times SL_N$ on its
homogeneous space $SL_N$
$$(g_1,g_2):\quad g \mapsto g_1gg_2^{-1},\qquad g,g_1,g_2 \in SL_N.$$
The action  by 'right shifts' is crucial in what follows, so we write $\xi
f$ instead of $(1 \otimes \xi)f$, $\xi \in U_q \frak{sl}_N$, $f \in{\Bbb
C}[SL_N]_q$. One has:
$$X_i^+t_{jk}=\left \{\begin{array}{ccl}t_{j \,k-1} &,& k=i+1 \\ 0 &,&{\rm
otherwise}\end{array}\right.,\qquad X_i^-t_{jk}=\left
\{\begin{array}{ccl}t_{j \,k+1} &,& k=i \\ 0 &,&{\rm
otherwise}\end{array}\right.,$$
$$H_it_{jk}=\left \{\begin{array}{ccl}t_{jk} &,& k=i \\ -t_{jk} &,& k=i+1 \\
0 &,&{\rm otherwise}\end{array}\right..$$

 The generators $t_{ij}$, $i,j=1,\ldots,N$, of ${\Bbb C}[SL_N]_q$ are just
the matrix elements of $\pi_1$. We are about to introduce the notation for
matrix elements of other fundamental representations of the quantum group
$SL_N$. Let $k$ be a natural number which does not exceed $N$, and consider
the representation $\pi_1^{\otimes k}$ of $U_q \frak{sl}_N$. Associate to
each collection $J$ of natural numbers $j_1<j_2<\ldots<j_k \le N$ the vector
$$u_J=u_{j_1}\wedge u_{j_2}\wedge \ldots \wedge u_{j_k}\stackrel{\rm
def}{=}\sum_{s \in S_k}(-q)^{l(s)}u_{j_{s(1)}}\otimes u_{j_{s(2)}}\otimes
\ldots \otimes u_{j_{s(k)}}.$$

 These vectors form a basis in the space of the representation
$\pi_1^{\wedge k}$ of $U_q \frak{sl}_N$. The matrix elements of
$\pi_1^{\wedge k}$ with respect to the basis $\{u_J \}$ are of the form
\begin{equation} t_{IJ}^{\wedge k}\stackrel{\rm def}{=}\sum_{s \in
S_k}(-q)^{l(s)}t_{i_1j_{s(1)}}\cdot t_{i_2j_{s(2)}}\cdot
\ldots \cdot t_{i_kj_{s(k)}},
\end{equation}
with $I=(i_1,i_2,\ldots,i_k)$, $J=(j_1,j_2,\ldots,j_k)$.

 Consider the element
\begin{equation} t \stackrel{\rm def}{=}t_{\{1,2,\dots,m
\}\{n+1,n+2,\ldots,N \}}^{\wedge m}.
\end{equation}

\medskip

\begin{lemma}
\begin{equation}\label{tijt}t_{ij}\cdot t=\left \{\begin{array}{ccl}qt \cdot
t_{ij} &,& i \le m \quad \&\quad j \le n \\ q^{-1}t \cdot t_{ij} &,& i>m
\quad \&\quad j>n \\ t \cdot t_{ij} &,& {\rm otherwise}\end{array}\right.,
\end{equation}
$$X_n^+t=t_{\{1,2,\dots,m \}\{n,n+2,\ldots,N \}}^{\wedge m},\qquad
X_n^-t=0,\qquad H_nt=-t.$$
\end{lemma}

\smallskip

 {\bf Proof.} The latter three equalities follow directly from the
definitions. The commutation relations (\ref{tijt}) is well known \cite{CP};
we present the proof for the reader's convenience in section 6. \hfill
$\Box$

\medskip

\begin{corollary} For any polynomial $f \in{\Bbb C}[t]$ one has:
\begin{equation}\label{Xf}X_j^+f(t)=\left
\{\begin{array}{ccl}t_{\{1,2,\dots,m \}\{n,n+2,\ldots,N \}}^{\wedge m}\cdot
\frac{\textstyle f(q^{-1}t)-f(t)}{\textstyle q^{-1}t-t} &,& j=n \\ 0 &,& j
\ne n \end{array}\right.,
\end{equation}
\begin{equation}\label{Hf}H_jf(t)=\left \{\begin{array}{ccl}-t{\textstyle
df(t)\over \textstyle dt} &,& j=n \\ 0 &,& j \ne n \end{array}\right.,\qquad
X_j^-f(t)=0,\;j=1,\ldots,N-1.
\end{equation}
\end{corollary}

\medskip

 Let ${\Bbb C}[SL_N]_{q,t}$ stand for the localization \index{localization}
of ${\Bbb C}[SL_N]_q$ with respect to the multiplicative system
$t,t^2,t^3,\ldots$.

 ${\Bbb C}[SL_N]_{q,t}$ has no zero divisors, its generators are $t^{-1}$,
$t_{ij}$, $i,j=1,\ldots,N$, and the relation list includes all the relations
which determine ${\Bbb C}[SL_N]_q$ and
\begin{equation}\label{tt^-1}{t^{-1}\cdot t_{\{1,2,\dots,m
\}\{n+1,n+2,\ldots,N \}}^{\wedge m}-1=0 \atop t_{\{1,2,\dots,m
\}\{n+1,n+2,\ldots,N \}}^{\wedge m}\cdot t^{-1}-1=0}
\end{equation}

 Apply the relations (\ref{Xf}), (\ref{Hf}) to equip ${\Bbb C}[SL_N]_{q,t}$
with a structure of covariant algebra in such a way that the canonical
embedding ${\Bbb C}[SL_N]_q \hookrightarrow{\Bbb C}[SL_N]_{q,t}$ becomes a
morphism of $U_q \frak{sl}_N$-modules:
$$X_j^+(t^{-1})=\left
\{\begin{array}{ccl}-q \cdot t_{\{1,2,\dots,m \}\{n,n+2,\ldots,N \}}^{\wedge
m}\cdot t^{-2} &,& j=n \\ 0 &,& j \ne n \end{array}\right.,$$
$$H_j(t^{-1})=\left \{\begin{array}{ccl}t^{-1} &,& j=n \\ 0 &,& j \ne n
\end{array}\right.,\qquad X_j^-(t^{-1})=0,\;j=1,\ldots,N-1.$$

 These rules are well defined, as one can easily see by applying $X_j^\pm$,
$H_j$, $j=1,\ldots,N-1$, to the left hand sides of (\ref{tt^-1}).

\medskip

 {\sc Remark 5.3.} One can also extend the structure of $U_q
\frak{sl}_N^{\rm op}$-module algebra in the same way. Thus, ${\Bbb
C}[SL_N]_{q,t}$ becomes a $U_q \frak{sl}_N^{\rm op}\otimes U_q
\frak{sl}_N$-module algebra.

\medskip \stepcounter{theorem}

 The following results of the present section are essentially due to M.
Noumi \cite{N}. They will be also refined in a subsequent section.

 Introduce the notation $J_{a \alpha}=\{n+1,n+2,\ldots,N \}\setminus
\{N+1-\alpha \}\cup \{a \}$ \index{$J_{a \alpha}$}.

\medskip

\begin{proposition}\label{emb} (cf. \cite{N}) The map $i:z_a^\alpha \mapsto
t^{-1}\cdot t_{\{1,2,\ldots,m \}J_{a \alpha}}$, $\alpha=1,\ldots,m$,
$a=1,\ldots,n$, admits a unique extension up to an embedding of $U_q
\frak{sl}_N$-module algebras $i:{\Bbb C}[{\rm Mat}_{mn}]_q
\hookrightarrow{\Bbb C}[SL_N]_{q,t}$.
\end{proposition}

\smallskip

 {\bf Proof.} We embed the algebras in question into the vector space
$(\widetilde{w}_0 \cdot U_q \frak{sl}_N)^*$, with $\widetilde{w}_0 \in{\Bbb
C}[SL_N]_q^*$ being the maximum length element of the quantum Weyl group
(see Appendix 1). More exactly, we restrict ourselves to the subspace ${\Bbb
F}\in(\widetilde{w}_0 \cdot U_q \frak{sl}_N)^*$ generated by $U_q
\frak{b}_-$-finite weight vectors
$$\{f \in(\widetilde{w}_0 \cdot U_q \frak{sl}_N)^*|\;{\rm dim}(U_q
\frak{b}_-f)<\infty,\quad H_if=\mu_if,\quad \mu_i \in{\Bbb Z},\quad
i=1,2,\ldots,N-1 \},$$
with $U_q \frak{b}_-$ being the standard Borel subalgebra of $U_q
\frak{sl}_N$. ${\Bbb F}$ is equipped with a structure of $U_q
\frak{sl}_N$-module algebra by the following identities derived from (A1.9):
$$\langle f_1f_2,\widetilde{w}_0 \xi \rangle=\langle
R_{\Bbb FF}\Delta(\xi)f_1 \otimes f_2,\widetilde{w}_0 \otimes
\widetilde{w}_0 \rangle,$$
$$\langle \xi f,\widetilde{w}_0 \eta \rangle=\langle f,\widetilde{w}_0 \eta
\xi \rangle \qquad \xi,\eta \in U_q \frak{sl}_N,\quad f,f_1,f_2 \in{\Bbb
F}.$$
(Here $R_{\Bbb FF}$ is the linear operator in ${\Bbb F \otimes F}$
determined by the universal R-matrix (see Appendix 1)).

 By a virtue of (A1.9), there exists an embedding of covariant algebras
${\Bbb C}[SL_N]_q \hookrightarrow{\Bbb F}$. It is worthwhile to note that the
element $t \in{\Bbb F}$ is invertible. (The proof of invertibility requires
some additional constructions. Given weight vector $f$ in a $U_q
\frak{sl}_N$-module ${\Bbb C}[SL_N]_q$, the sequence $c_k(f)=\langle
ft^k,\widetilde{w}_0 \rangle$ satisfies a difference equation of order one
derived from (A1.9):
$$\langle ft^k,\widetilde{w}_0 \rangle=\langle R_{{\Bbb C}[SL_N]_q \,{\Bbb
C}[SL_N]_q}(ft^{k-1}\otimes t),\widetilde{w}_0 \otimes \widetilde{w}_0
\rangle,\qquad k \in{\Bbb N}.$$
That is, by a virtue of (A1.6),
$$\langle ft^k,\widetilde{w}_0 \rangle=\langle q^{-{H_0 \otimes H_0 \over
(H_0,H_0)}}(ft^{k-1}\otimes t),\widetilde{w}_0 \otimes \widetilde{w}_0
\rangle,\qquad k \in{\Bbb N}.$$
This allows one to extend the linear functional $\widetilde{w}_0$ from
${\Bbb C}[SL_N]_q$ onto ${\Bbb C}[SL_N]_{q,t}$. We shall use in the sequel
exactly this extension, together with the following pairing of ${\Bbb
C}[SL_N]_{q,t}$ and $\widetilde{w}_0U_q \frak{sl}_N$:
$$\langle f,\widetilde{w}_0 \xi \rangle \stackrel{\rm def}{=}\langle \xi
f,\widetilde{w}_0 \rangle,\qquad \xi \in U_q \frak{sl}_N,\quad f \in{\Bbb
C}[SL_N]_{q,t}.$$
Its non-degeneracy follows from the fact that its restriction onto ${\Bbb
C}[SL_N]_q \times \widetilde{w}_0U_q \frak{sl}_N$ is non-degenerate. In fact,
if $\langle f,\widetilde{w}_0 \xi \rangle=0$ for all $\xi \in U_q
\frak{sl}_N$, then $\langle f \cdot t^j,\widetilde{w}_0 \xi \rangle=0$ since
$\langle f \cdot t^j,\widetilde{w}_0 \xi \rangle=\langle R_{{\Bbb
C}[SL_N]_{q,t}\,{\Bbb C}[SL_N]_{q,t}}\Delta(\xi)(f \otimes
t^j),\widetilde{w}_0 \otimes \widetilde{w}_0 \rangle$.  Thus we get an
embedding ${\Bbb C}[SL_N]_{q,t}\hookrightarrow(\widetilde{w}_0U_q
\frak{sl}_N)^*$. What remains is to note that the image of $t^{-1}$ under
this embedding is in ${\Bbb F}$.)

 Consider the onto linear map $j:\widetilde{w}_0U_q \frak{sl}_N \to V_-(0)$,
$j:\widetilde{w}_0 \xi \mapsto S(\xi)v_-(0)$, $\xi \in U_q \frak{sl}_N$,
with $v_-(0)$ being the generator of the $U_q \frak{sl}_N$-module $V_-(0)$
dual to ${\Bbb C}[{\rm Mat}_{mn}]_q$ (see \cite{SV2}). It is easy to prove
that the adjoint linear map $j^*:{\Bbb C}[{\rm Mat}_{mn}]_q
\hookrightarrow{\Bbb F}$ is an embedding of $U_q \frak{sl}_N$-module
algebras. Let us agree not to distinguish between the $U_q
\frak{sl}_N$-module algebras ${\Bbb C}[{\rm Mat}_{mn}]_q$, ${\Bbb
C}[SL_N]_{q,t}$ and their images under the above embedding into ${\Bbb F}$.
In view of propositions \ref{zaalph}, \ref{gen}, it suffices to prove that
$t^{-1}t_{\{1,2,\ldots,m \}J_{a \alpha}}\in{\Bbb C}[{\rm Mat}_{mn}]_q$,
$\alpha=1,2,\ldots,m$, $a=1,2,\ldots,n$. For that, we need only to establish
that $t^{-1}t_{\{1,2,\ldots,m \}J_{a \alpha}}$ are orthogonal to the kernel
of $j$ with respect to the above pairing. This kind of orthogonality follows
from
$$((K_j^{\pm 1}-1)\otimes 1)t^{-1}t_{\{1,2,\ldots,m \}J_{a
\alpha}}=(F_j \otimes 1)t^{-1}t_{\{1,2,\ldots,m \}J_{a \alpha}}=0,$$
$$(E_i \otimes 1)t^{-1}t_{\{1,2,\ldots,m \}J_{a \alpha}}=0,\qquad
j=1,2,\ldots,N-1,\quad i=1,2,\ldots,n-1,n+1,\ldots,N,$$
in view of the definitions of the $U_q \frak{sl}_N$-module $V_-(0)$ and
$\widetilde{w}_0$. \hfill $\Box$

\medskip

 Let us agree not to distinguish the elements of ${\Bbb C}[{\rm
Mat}_{mn}]_q$ and their images under the embedding $i:{\Bbb C}[{\rm
Mat}_{mn}]_q \hookrightarrow{\Bbb C}[SL_N]_{q,t}$.

\medskip

\begin{lemma} For all $1 \le a<b \le n$, $1 \le \alpha<\beta \le m$,
$$t^{-1}\cdot t_{\{1,2,\ldots,m
\}\{a,b,\ldots,\widehat{N+1-\beta},\ldots,\widehat{N+1-\alpha},\ldots,N
\}}=z_a^\alpha z_b^\beta-qz_a^\beta z_b^\alpha.$$
\end{lemma}

\smallskip

 {\bf Proof.} In the same way as in the proof of proposition \ref{emb}, one
can establish that $t^{-1}\cdot t_{\{1,2,\ldots,m
\}\{a,b,\ldots,\widehat{N+1-\beta},\ldots,\widehat{N+1-\alpha},\ldots,N
\}}\in{\Bbb C}[{\rm Mat}_{mn}]_q \subset{\Bbb F}$. What remains is to
express this element in terms of the generators of ${\Bbb C}[{\rm
Mat}_{mn}]_q$. It is a weight vector of the $U_q \frak{sl}_N$-module ${\Bbb
C}[{\rm Mat}_{mn}]_q$. A computation of the weight yields
$$t^{-1}\cdot t_{\{1,2,\ldots,m
\}\{a,b,\ldots,\widehat{N+1-\beta},\ldots,\widehat{N+1-\alpha},\ldots,N
\}}=c_1z_a^\beta z_b^\alpha+c_2z_a^\alpha z_b^\beta,$$
with $c_1$, $c_2$ being the elements of the ground field ${\Bbb
C}(q^{1/s})$. When computing the constants $c_1$, $c_2$, one can restrict
oneself to the special case $m=n=2$ by passing from the algebra ${\Bbb
C}[SL_N]_{q,t}$ to the corresponding factor algebra. In the special case
$m=n=2$ the result in question is accessible via a direct calculation
\cite{N}.  \hfill $\Box$

\medskip

\begin{corollary}\label{Xnzaa}
$$X_n^+z_a^\alpha=\left \{\begin{array}{ccl}-q^{-1/2}z_a^mz_n^\alpha &,& a
\ne n \quad \&\quad \alpha \ne m \\ -q^{-1/2}(z_n^m)^2 &,& a=n \quad \&\quad
\alpha=m \\ -z_n^mz_a^\alpha &,& {\rm otherwise}\end{array}\right.$$
The latter relation, together with  proposition 2.1 and the relations
$$X_n^-z_a^\alpha=\left \{\begin{array}{ccl}q^{1/2} &,& a=n \quad \&\quad
\alpha=m \\ 0 &,& {\rm
otherwise}\end{array}\right.,$$
$$H_nz_a^\alpha=\left \{\begin{array}{ccl}2z_a^\alpha &,& a=n \quad \&\quad
\alpha=m \\ z_a^\alpha &,& a=n \quad \&\quad \alpha \ne m \quad{\rm or}\quad
a \ne n \quad \&\quad \alpha=m \\ 0 &,& {\rm otherwise}\end{array}\right.$$
describe the action of $U_q \frak{sl}_N$ in ${\Bbb C}[{\rm Mat}_{mn}]_q$.
\end{corollary}

\bigskip

\section{The quantum principal homogeneous space}

 In the case $q=1$ the matrix ball ${\Bbb U}$ as a homogeneous space of the
group $SU_{mn}$ is isomorphic to $S(U_n \times U_m)\setminus SU_{mn}$. A
straightforward generalization of this statement is derived via a
replacement of $S(U_n \times U_m)\setminus SU_{mn}$ by $S(U_n \times
U_m)\setminus \widetilde{X}$, with $\widetilde{X}$ being some principal
homogeneous space of the group $SU_{mn}$. We construct a quantum principal
homogeneous space in such a way that the isomorphism ${\Bbb U}\simeq S(U_n
\times U_m)\setminus \widetilde{X}$ is valid in the quantum case.

 Consider the element $\widetilde{w}_0 \in{\Bbb C}[SL_N]_q^*$ of the quantum
Weyl group (see Appendix 1). It follows from the invertibility of this
element that the pairing of ${\Bbb C}[SL_N]_q$ and $\widetilde{w}_0U_q
\frak{sl}_n$ is non-degenerate, and hence there exists a unique antilinear
operator $*$ in ${\Bbb C}[SL_N]_q$ such that
\begin{equation} \langle f^*,\widetilde{w}_0 \xi \rangle=\overline{\langle
f,\widetilde{w}_0(S(\xi))^*\rangle}
\end{equation}
for all $f \in{\Bbb C}[SL_N]_q$, $\xi \in U_q \frak{su}_{mn}$.

\medskip

\begin{proposition} The map $*$ is an antilinear involution:
\begin{equation}\label{*}f^{**}=f,\qquad (f_1f_2)^*=f_2^*f_1^*,\qquad
f,f_1,f_2 \in{\Bbb C}[SL_N]_q.
\end{equation}
\end{proposition}

\smallskip

 {\bf Proof.} It follows from the well known properties of a universal
R-matrix (see \cite{CP}) and the definition of involution $*$ that
\begin{equation}\label{R*S} R^{*\otimes*}=R_{21},\qquad S \otimes S(R)=R,
\end{equation}
with $R_{21}$ is derivable from $R$ via a displacement of tensor multiples.
The second one of the identities (\ref{*}) follows from the relations (A1.9)
and (\ref{R*S}), and the first one follows from $(S((S \xi)^*))^*=\xi$,
which is valid for all $\xi \in U_q \frak{su}_{mn}$. \hfill $\Box$

\medskip

 Evidently, (\ref{*H}) is valid for the $*$-algebra ${\rm
Pol}(\widetilde{X})_q=({\Bbb C}[SL_N]_q,*)$ \index{${\rm
Pol}(\widetilde{X})_q$}. That is, ${\rm Pol}(\widetilde{X})_q$ is a
covariant $*$-algebra. We call it the polynomial algebra on the quantum
principal homogeneous space \index{polynomial algebra on the quantum
principal homogeneous space}. (In the classical case $q \to 1$ one has
$\widetilde{X}=\widetilde{w}_0 SU_{mn}\subset SL_N$.) We are to produce an
explicit formula for $t_{ij}^*$, $i,j=1,\ldots,N$. Let us consider the
covector representation $\pi _{N-1}$ of $U_q \frak{sl}_{N}$ defined in the
base $v^1$,$v^2$,...$v^N$ by the formulae from section 2. In the next lemma
we express all the matrix elements $t_{ij}'\in{\Bbb C}[SL_N]_q$ of $\pi
_{N-1}$ $$\pi_{N-1}(\xi)v^j=\sum_i{\langle t_{ij}',\xi \rangle v^i}$$
 in terms of generators $t_{ij}$.

\medskip

\begin{lemma}\label{tij'd} For all $i,j=1,\ldots,N$,
$$t_{ij}'={\rm det}_qT_{i'j'},\qquad with \quad i'=N+1-i,\;j'=N+1-j.$$
\end{lemma}

\smallskip

 {\bf Proof.} Let $L$ be the linear span of $t_{ij}'\in{\Bbb C}[SL_N]_q$.
Evidently, $L$ is a {\sl simple} $U_q \frak{sl}_N^{\rm op}\otimes U_q
\frak{sl}_N$-module, and the map $t_{ij}'\mapsto{\rm det}_qT_{i'j'}$,
$i,j=1,\ldots,N$, admits an extension up to an endomorphism $\varphi$ of
this simple module. Hence, $\varphi={\tt const}\cdot 1$. On the other hand,
$\langle t_{ij},1 \rangle=\langle t_{ij}',1 \rangle=\delta_{ij}$. Thus we
have $\langle t_{NN}',1 \rangle=\langle{\rm det}_qT_{N'N'},1 \rangle=1$, and
so $\varphi=1$ and $t_{ij}'={\rm det}_qT_{i'j'}$ for all $i,j=1,\ldots,N$.
\hfill $\Box$

\medskip

 Compare the matrices of the operators $\pi_1(w_0)$ and $\pi_{N-1}(w_0)$.

\medskip

\begin{lemma}\label{tij'} For all $i,j=1,\ldots,N$,
\begin{equation} \langle t_{ij},\widetilde{w}_0 \rangle=\langle
t_{ij}',\widetilde{w}_0 \rangle={\tt const}\cdot(-q)^{-i}\cdot
\delta_{i+j,N+1},
\end{equation}
with {\tt const} being an element of the ground field independent of $i$,
$j$.
\end{lemma}

\smallskip

 {\bf Proof.} In the case $N=2$ the desired statement is well known
\cite{SoV2, KR}. The general case is reducible to this one via (A1.10),
(A1.11) and the following reduced decompositions of the full permutation
$w_0=(N,N-1,\ldots,2,1)\in S_N$:
$$w_0=s_{N-1}(s_{N-2}s_{N-1})\ldots(s_1s_2 \ldots
s_{N-1})=s_1(s_2s_1)\ldots(s_{N-1}s_{N-2}\ldots s_1).$$
(Note that for any orthogonal basis $\{I_k \}_{k=1}^{N-1}$ in the Cartan
subalgebra one has
$$\pi_1 \left(\sum_k{I_k^2 \over(I_k,I_k)}\right)={\tt const}'\cdot I,\qquad
\pi_{N-1}\left(\sum_k{I_k^2 \over(I_k,I_k)}\right)={\tt const}''\cdot I,$$
with $I$ being the identity operator. This is because the element
$\displaystyle \sum_k{I_k^2 \over(I_k,I_k)}$ is an invariant of the Weyl
group action. That is, by a slight misuse of the notation,
$$\langle t_{ij},\sum_k{I_k^2 \over(I_k,I_k)}\rangle={\tt const}'\cdot
\delta_{ij},\qquad \langle t_{ij}',\sum_k{I_k^2 \over(I_k,I_k)}\rangle={\tt
const}''\cdot \delta_{ij}.$$
In this setting ${\tt const}'={\tt const}''$ since $\langle
t_{ij}',\xi \rangle=\langle t_{ij},\omega(\xi)\rangle$, $\xi \in U_q
\frak{sl}_N$, $i,j=1,\ldots,N$, with $\omega:U_q \frak{sl}_N \to U_q
\frak{sl}_N$ being the automorphism given by $\omega(K_j^{\pm
1})=K_{N-j}^{\pm 1}$, $\omega(E_j)=E_{N-j}$, $\omega(F_j)=F_{N-j}$.) \hfill
$\Box$

\medskip

\begin{lemma}\label{l##}\hfill \\
1) There exists a unique antilinear involution $\#$ in ${\Bbb C}[SL_N]_q$
such that for all $f \in{\Bbb C}[SL_N]_q$, $\xi \in U_q \frak{su}_{nm}$,
$$\langle f^\#,\xi \rangle=\overline{\langle f,(S(\xi))^*\rangle}.$$
2) For all $i,j=1,\ldots,N$,
\begin{equation}\label{##}t_{ij}^\#={\rm sign}\left((n-i+{\scriptstyle 1
\over \scriptstyle 2})(n-j+{\scriptstyle 1 \over \scriptstyle
2})\right)(-q)^{j-i}{\rm det}_qT_{ij}.
\end{equation}
\end{lemma}

\smallskip

 {\bf Proof.} The desired statement is a well known fact in quantum group
theory. It follows from the results of \cite{RTF} where the involution
(\ref{##}) was initially considered. (It was also noted in this paper that
the $*$-algebra ${\Bbb C}[SU_{nm}]_q \stackrel{\rm def}{=}({\Bbb
C}[SL_N]_q,\#)$ is a Hopf $*$-algebra.) \hfill $\Box$

\medskip

\begin{proposition} For all $i,j=1,\ldots,N$,
\begin{equation}\label{*Tij}t_{ij}^*={\rm sign}\left((i-m-{\scriptstyle 1
\over \scriptstyle 2})(n-j+{\scriptstyle 1 \over \scriptstyle
2})\right)(-q)^{j-i}{\rm det}_qT_{ij}.
\end{equation}
\end{proposition}

\smallskip

 {\bf Proof.} The pairing considered above is non-degenerate and allows one
to embed $U_q \frak{sl}_N$ into ${\Bbb C}[SL_N]_q^*$. Let $L$ be the
antirepresentation of ${\Bbb C}[SL_N]_q^*$ in the space ${\Bbb C}[SL_N]_q$
given by
$$\langle L(\xi)f,\eta \rangle=\langle f,\xi \eta \rangle,\qquad f \in{\Bbb
C}[SL_N]_q,\;\xi,\eta \in{\Bbb C}[SL_N]_q^*.$$
Now compare the definitions for involutions $*$ and $\#$ to get
$$t_{ij}^*=L(\widetilde{w}_0^{-1})\cdot(L(\widetilde{w}_0)t_{ij})^\#,\qquad
i,j=1,\ldots,N.$$

 On the other hand, it follows from lemma \ref{tij'} that
$$L(\widetilde{w}_0)t_{ij}={\tt const}\cdot(-q^{-1})^i \cdot t_{N+1-i,j},$$
$$L(\widetilde{w}_0^{-1})t_{ij}'={1 \over{\tt const}}\cdot(-q)^{N-i+1}\cdot
t_{N+1-i,j}',$$
and, by a virtue of lemmas \ref{tij'd}, \ref{l##},
$$t_{ij}^\#={\rm sign}\left((n-i+{\scriptstyle 1 \over \scriptstyle
2})(n-j+{\scriptstyle 1 \over \scriptstyle
2})\right)(-q)^{j-i}t_{N+1-i,N+1-j}'.$$
Hence,
$$t_{ij}^*={\tt
const}\cdot(-q^{-1})^iL(\widetilde{w}_0^{-1})t_{N+1-i,j}^\#=$$
$$={\tt const}\cdot(-q^{-1})^iL(\widetilde{w}_0^{-1})\cdot{\rm
sign}\left((i-m-{\scriptstyle 1 \over \scriptstyle 2})(n-j+{\scriptstyle 1
\over \scriptstyle 2})\right)(-q)^{i+j-N-1}\cdot
t_{i,N+1-j}'=$$
$$=(-q^{-1})^i \cdot{\rm sign}\left((i-m-{\scriptstyle 1 \over \scriptstyle
2})(n-j+{\scriptstyle 1 \over \scriptstyle
2})\right)(-q)^{i+j-N-1}\cdot(-q)^{N-i+1}\cdot t_{N+1-i,N+1-j}'.$$
What remains is to apply lemma \ref{tij'd}. \hfill $\Box$

\medskip

 We use in the sequel the results of Ya. Soibelman concerning the quantum
group $SU_N$ \cite{So1, SoV3}. Remind the definitions.

 Equip the Hopf algebras $U_q \frak{sl}_N$ and ${\Bbb C}[SL_N]_q$ with
antilinear involutions given by
$$(K_j^{\pm 1})^\star=K_j^{\pm 1},\qquad E_j^\star=K_jF_j,\qquad
F_j^\star=E_jK_j^{-1},\qquad j=1,\ldots,N-1,$$
$$\langle f^\star,\xi \rangle=\overline{\langle f,(S(\xi))^\star
\rangle},\qquad f \in{\Bbb C}[SL_N]_q,\;\xi \in U_q \frak{sl}_N.$$
The associated Hopf $*$-algebras are denoted by
$$U_q \frak{su}_N \stackrel{\rm def}{=}(U_q \frak{sl}_N,\star),\qquad {\Bbb
C}[SU_N]_q=({\Bbb C}[SL_N]_q,\star).$$
\index{$U_q \frak{su}_N$} \index{${\Bbb C}[SU_N]_q$} Similarly to (\ref{##}),
(\ref{*Tij}), one has
\begin{equation}\label{tijst}t_{ij}^\star=(-q)^{j-i}{\rm det}_qT_{ij},
\end{equation}
together with the following

\medskip

\begin{lemma}\label{twedgest} For all $k=1,\ldots,N-1$,
$$\left(t_{\{1,\ldots,k \}\{N-k+1,\ldots,N \}}^{\wedge k}
\right)^\star=(-q)^{k(N-k)}\cdot t_{\{k+1,\ldots,N \}\{1,\ldots,N-k
\}}^{\wedge(N-k)}.$$
\end{lemma}

\medskip

Introduce the notation $t=t_{\{1,2,\ldots,m \}\{n+1,n+2,\ldots,N
\}}^{\wedge m}$, $x=(-q)^{mn}\cdot t_{\{1,2,\ldots,m \}\{n+1,n+2,\ldots,N
\}}^{\wedge m}\cdot t_{\{m+1,m+2,\ldots,N \}\{1,2,\ldots,n \}}^{\wedge n}$
for the elements of a crucial importance in the function theory in quantum
matrix ball. Note that in view of (\ref{*Tij}) one has $x=tt^*$.

\medskip

\begin{lemma}\label{tcr}
\begin{equation}\label{tijt_}t_{ij}\cdot t=\left \{\begin{array}{ccl}qt \cdot
t_{ij} &,& i \le m \quad \&\quad j \le n \\ q^{-1}t \cdot t_{ij} &,& i>m
\quad \&\quad j>n \\ t \cdot t_{ij} &,& {\rm otherwise}\end{array}\right.,
\end{equation}
\begin{equation}\label{tijt*}t_{ij}\cdot t^*=\left
\{\begin{array}{ccl}qt^*\cdot t_{ij} &,& i \le m \quad \&\quad j \le n \\
q^{-1}t^*\cdot t_{ij} &,& i>m \quad \&\quad j>n \\ t^*\cdot t_{ij} &,& {\rm
otherwise}\end{array}\right..
\end{equation}
\end{lemma}

\medskip

 {\bf Proof.} (\ref{tijt_}) can be easily verified in the special case $i \in
\{m,m+1 \}$, $j \in \{n,n+1 \}$. The biinvariance of $t$
$$(1 \otimes E_i)t=(1 \otimes F_i)t=(E_j \otimes 1)t=(F_j \otimes
1)t=0,\qquad i \ne n,\;j \ne m$$
allows one to reduce the general case to the above special case via an
application of the operators
$$1 \otimes E_i,\quad 1 \otimes F_i,\quad i \ne n;\qquad E_j \otimes
1,\quad F_j \otimes 1,\quad j \ne m$$
to each side of (\ref{tijt_}). A similar argument  proves also (\ref{tijt*}).
\hfill $\Box$

\medskip

 Now lemmas \ref{twedgest}, \ref{tcr} imply

\medskip

\begin{corollary} $x=tt^*=t^*t$, and for every polynomial $f \in{\Bbb
C}[x]$,
\begin{equation}\label{tijfx}t_{ij} \cdot f(x)=\left \{\begin{array}{ccl}
f(q^2x)t_{ij} &,& i \le m \quad \&\quad j \le n \\ f(q^{-2}x)t_{ij} &,& i>m
\quad \&\quad j>n \\ f(x)t_{ij} &,& {\rm otherwise}\end{array}\right..
\end{equation}
\end{corollary}

\medskip

\begin{proposition} For every polynomial $f \in{\Bbb C}[x]$,
\begin{equation}\label{Xn+f}X_n^+f(x)=(X_n^+x)\cdot
\frac{f(q^{-2}x)-f(x)}{q^{-2}x-x},
\end{equation}
\begin{equation}\label{Xn-f}
X_n^-f(x)=\frac{f(q^{-2}x)-f(x)}{q^{-2}x-x}\cdot(X_n^-x).
\end{equation}
\end{proposition}

\smallskip

 {\bf Proof.} It follows from the explicit formulae which define the action
of the operators $X_n^\pm$, $H_n$ on the genrators $t_{ij}$ of ${\Bbb
C}[SL_N]_q$ and the covariance of the latter algebra that
\begin{equation}\label{Xn+x}X_n^+x=q^{-1/2}(-q)^{mn}\cdot t_{\{1,2,\ldots,m
\}\{n,n+2,\ldots,N \}}^{\wedge m}\cdot t_{\{m+1,m+2,\ldots,N
\}\{1,2,\ldots,n \}}^{\wedge n},
\end{equation}
\begin{equation}\label{Xn-x}X_n^-x=q^{-1/2}(-q)^{mn}\cdot t_{\{1,2,\ldots,m
\}\{n+1,n+2,\ldots,N \}}^{\wedge m}\cdot t_{\{m+1,m+2,\ldots,N
\}\{1,2,\ldots,n-1,n+1 \}}^{\wedge n}.
\end{equation}

 It follows from (\ref{tijfx}), (\ref{Xn+x}), (\ref{Xn-x}) that
$$(X_n^+x)\cdot x=q^2x \cdot(X_n^+x),\qquad(X_n^-x)\cdot x=q^{-2}x
\cdot(X_n^-x).$$
Hence, (\ref{Xn+f}), (\ref{Xn-f}) are valid for all monomials $f=x^k$, $k
\in{\Bbb Z}_+$. \hfill $\Box$

\medskip

 Let ${\rm Pol}(\widetilde{X})_{q,x}$ stand for a localization of the
integral domain ${\rm Pol}(\widetilde{X})_q$ with respect to the
multiplicative system $x,x^2,x^3,\ldots$. An involution in ${\rm
Pol}(\widetilde{X})_{q,x}$ is imposed in a natural way:
$(x^{-1})^*=x^{-1}$.

 Apply (\ref{Xn+f}), (\ref{Xn-f}) to equip ${\rm Pol}(\widetilde{X})_{q,x}$
with a structure of $U_q \frak{su}_{nm}$-module algebra in such a way that
the canonical embedding ${\rm Pol}(\widetilde{X})_q \hookrightarrow{\rm
Pol}(\widetilde{X})_{q,x}$ becomes a morphism of $U_q
\frak{su}_{nm}$-modules:
$$H_j(x^{-1})=0,\qquad X_j^+(x^{-1})=\left
\{\begin{array}{ccl}-q^2(X_j^+x)x^{-2} &,& j=n \\ 0 &,& j \ne n
\end{array}\right.,$$
$$X_j^-(x^{-1})=\left \{\begin{array}{ccl}-q^2x^{-2}(X_j^-x) &,& j=n \\ 0
&,& j \ne n \end{array}\right..$$
This structure of $U_q \frak{sl}_N$-module algebra is well defined, as
one can easily verify just as for a similar statement in the previous
section. The relation $(\xi f)^*=(S(\xi))^*f^*$, $\xi \in U_q \frak{sl}_N$,
is valid both for $f=x^{-1}$ and $f \in{\rm Pol}(\widetilde{X})_q$. Hence it
is valid for all $f \in{\rm Pol}(\widetilde{X})_{q,x}$.

 Just as in remark 5.3, note that ${\rm Pol}(\widetilde{X})_{q,x}$ is a $U_q
\frak{sl}_N^{\rm op}\otimes U_q \frak{sl}_N$-module algebra.

\medskip

\begin{proposition} The map
\begin{equation}\label{cI}{\cal I}:z_a^\alpha \mapsto
t^{-1}t_{\{1,2,\ldots,m \}J_{a \alpha}}
\end{equation}
with $J_{a \alpha}=\{n+1,n+2,\ldots,N \}\setminus \{N+1-\alpha \}\cup \{a
\}$, is uniquely extendable up to an embedding  of $U_q
\frak{su}_{nm}$-module *-algebras ${\cal I}:{\rm Pol}({\rm
Mat}_{mn})_q \hookrightarrow{\rm Pol}(\widetilde{X})_{q,x}$.
\end{proposition}

\smallskip

 {\bf Proof.} The uniqueness of the embedding ${\cal I}$ follows from
proposition \ref{gen}. Prove its existence. Consider the embedding of $U_q
\frak{sl}_{N}$-module algebras $i:{\Bbb C}[{\rm Mat}_{mn}]_q
\hookrightarrow{\Bbb C}[SL_N]_{q,t}$ (see section 5) and a similar embedding
$\overline{i}:{\Bbb C}[\overline{\rm Mat}_{mn}]_q \hookrightarrow{\Bbb
C}[SL_N]_{q,t^*}$
$$\overline{i}f=(if^*)^*,\qquad f \in{\Bbb C}[\overline{\rm Mat}_{mn}]_q.$$
(We use the embeddings ${\Bbb C}[\overline{\rm Mat}_{mn}]_q \subset{\rm
Pol}({\rm Mat}_{mn})_q$, ${\Bbb C}[SL_N]_{q,t}\subset{\rm
Pol}(\widetilde{X})_{q,x}$.)

 Let ${\cal I}$ be such a linear operator ${\cal I}:{\rm Pol}({\rm
Mat}_{mn})_q \to{\rm Pol}(\widetilde{X})_{q,x}$ that ${\cal I}:f_- \cdot f_+
\mapsto \overline{i}(f_-) \cdot i(f_+)$, $f_-\in{\Bbb C}[\overline{\rm
Mat}_{mn}]_q$, $f_+\in{\Bbb C}[{\rm Mat}_{mn}]_q$. By our construction,
${\cal I}$ is a morphism of $U_q \frak{su}_{nm}$-modules and satisfies
(\ref{cI}). Prove that it is a homomorphism of $*$-algebras. It suffices to
show that (in the notation of section 4)
$$({\cal I}z_b^\beta)^*({\cal I}z_a^\alpha)=mPR_{{\cal I}{\Bbb
C}[\overline{\rm Mat}_{mn}]_q \:{\cal I}{\Bbb C}[{\rm Mat}_{mn}]_q}(({\cal
I}z_b^\beta)^*\otimes({\cal I}z_a^\alpha))$$
for all $a,b=1,2,\ldots,n$; $\alpha,\beta=1,2,\ldots,m$. For that, it
suffices to establish
$${\cal I}((z_b^\beta)^*)t^{*j}\cdot t^k{\cal I}(z_a^\alpha)=q^{{\tt
const}\cdot j \cdot k}mPR_{{\Bbb C}[SL_N]_{q,t^*}{\Bbb C}[SL_N]_{q,t}}({\cal
I}((z_b^\beta)^*)t^{*j}\otimes t^k{\cal I}(z_a^\alpha)),$$
with $j,k \in{\Bbb Z}$, and {\tt const} being determined by the equation
$H_0 \otimes H_0(t \otimes t^*)={\tt const}\cdot(H_0,H_0)\cdot t \otimes
t^*$. We may restrict ourselves to the special case $j,k \in{\Bbb N}$ since
$$t^{-k}(q^{{\tt const}\cdot j \cdot k}mPR_{{\Bbb C}[SL_N]_{q,t^*}{\Bbb
C}[SL_N]_{q,t}}({\cal I}((z_b^\beta)^*)t^{*j}\otimes t^k{\cal
I}(z_a^\alpha)))t^{*-j}$$
is a Laurent polynomial of $q^{k/s}$, $q^{j/s}$. In the above special case
${\cal I}((z_b^\beta)^*)t^{*j}$, $ t^k{\cal I}(z_a^\alpha)$,
$a,b=1,2,\ldots,n$, $\alpha,\beta=1,2,\ldots,m$, are the matrix elements of
finite dimensional representations of $U_q \frak{sl}_N$. What remains is to
apply the well known \cite{CP} R-matrix commutation relations between those
matrix elements, (A1.6), and the relations $(F_i \otimes 1)(t^k{\cal
I}(z_a^\alpha))=(E_i \otimes 1)({\cal I}((z_b^\beta)^*) t^{*j})=0$,
$a,b=1,2,\ldots,n$, $\alpha,\beta=1,2,\ldots,m, i \ne m$.

 So we need only to prove the injectivity of the homomorphism ${\cal I}$ via
passage to the 'improper specialization $q=1$' and observing that the
corresponding homomorphism in the case $q=1$ is injective (see the proof of
proposition \ref{gen} where a similar well known argument was used). \hfill
$\Box$

\medskip

 To conclude, extend the embedding ${\Bbb C}[SL_N]_q
\hookrightarrow(\widetilde{w}_0U_q \frak{sl}_N)^*$ initially constructed in
the proof of proposition \ref{emb}, up to an embedding ${\rm
Pol}(\widetilde{X})_{q,x}\hookrightarrow(\widetilde{w}_0U_q \frak{sl}_N)^*$.
One can verify in the same way as in section 5 that for every weight vector
$f \in{\rm Pol}(\widetilde{X})_{q,x}$, the numbers $c_{jk}=\langle
t^{*j}ft^k,\widetilde{w}_0 \rangle$, $j,k \in{\Bbb Z}_+$, satisfy the system
of order one difference equations
$$\left \{\begin{array}{ccl}\langle t^{*j}ft^k,\widetilde{w}_0 \rangle &=&
\langle q^{-{H_0 \otimes H_) \over(H_0,H_0)}}t^* \otimes
t^{*(j-1)}ft^k,\widetilde{w}_0 \otimes \widetilde{w}_0 \rangle \\ \langle
t^{*j}ft^k,\widetilde{w}_0 \rangle &=& \langle q^{-{H_0 \otimes H_)
\over(H_0,H_0)}}t^{*j}ft^{k-1}\otimes t,\widetilde{w}_0 \otimes
\widetilde{w}_0 \rangle \end{array}\right..$$
Extend the linear functional $\widetilde{w}_0$ onto ${\rm
Pol}(\widetilde{X})_{q,x}$ via the above system of difference equations and
set up $\langle f,\widetilde{w}_0 \xi \rangle=\langle \xi f,\widetilde{w}_0
\rangle$, $\xi \in U_q \frak{sl}_N$, $f \in{\rm Pol}(\widetilde{X})_{q,x}$
(see the proof of proposition \ref{emb}).

 A non-degenerate pairing supplies an embedding ${\rm
Pol}(\widetilde{X})_{q,x}\hookrightarrow(\widetilde{w}_0U_q \frak{sl}_N)^*$.
According to an approach of V. Drinfeld, we use the term 'function on a
quantum group' to denote the linear functions on $U_q \frak{sl}_N$ (more
exactly, on $U_h \frak{sl}_N$, see \cite{D}). The basic distinction of our
approach to construction of algebras of functions on quantum homogeneous
spaces is in replacing the space $U_q \frak{sl}_N$ with $\widetilde{w}_0 U_q
\frak{sl}_N$.

\bigskip

\section{Algebras of finite functions}

 From now on we assume that $0<q<1$ and use ${\Bbb C}$ as a ground field. We
also keep the above notation for the Hopf algebras and covariant algebras
involved, together with their descriptions in terms of generators and
relations.

 Among the principal tools in harmonic analysis, one should mention the
algebra of finite functions and the $L^2$ space being its completion. Now
turn to a construction of covariant $*$-algebra $D(\widetilde{X})_q$ of
finite functions \index{finite functions} on the quantum principal
homogeneous space $\widetilde{X}$.

 We refer to Appendix 2 for a construction of a $*$-representation
$\widetilde{\Pi}$ of ${\rm Pol}(\widetilde{X})_q$ with
$\widetilde{\Pi}(x)\ne 0$ and ${\tt spec} \widetilde{\Pi}(x)=q^{-2{\Bbb
Z}_+}$.  For a function $f$ on $q^{-2{\Bbb Z}_+}$ with a finite support, a
bounded linear operator $f(\widetilde{\Pi}(x))$ is well defined. Our
immediate intention is to add the elements of the form $f(x)$, ${\tt supp}
f(x)\subset q^{-2{\Bbb Z}_+}$, to ${\rm Pol}(\widetilde{X})_q$.

 Consider the algebra of polynomial functions and the algebra of functions
with finite support inside $q^{2{\Bbb Z}}$; let ${\Bbb F}$ stand for
their sum. Note that ${\Bbb F}$ admits the involution $f(x)\mapsto
\overline{f(x)}$.

 Consider the ${\Bbb C}[x]$-bimodules ${\rm Pol}(\widetilde{X})_q$, ${\Bbb
F}$, and form their tensor product ${\cal F}={\rm Pol}(\widetilde{X})_q
\otimes_{{\Bbb C}[x]}{\Bbb F}$. One can deduce from (\ref{tijfx}) that there
exists a canonical isomorphism ${\cal F}\simeq{\Bbb F}\otimes_{{\Bbb
C}[x]}{\rm Pol}(\widetilde{X})_q$. This isomorphism, together with the
multiplication laws
$${\rm Pol}(\widetilde{X})_q^{\otimes 2}\to {\rm
Pol}(\widetilde{X})_q,\qquad{\Bbb F}\otimes{\Bbb F}\to{\Bbb F},$$
is used to equip ${\cal F}$ with such structure of $*$-algebra that the
embeddings $${\rm Pol}(\widetilde{X})_q \hookrightarrow{\cal F},\qquad{\Bbb
F}\hookrightarrow{\cal F}.$$ are  homomorphisms of algebras.

\medskip

\begin{proposition} There exists a unique extension of the structure of
covariant $*$-algebra from ${\rm Pol}(\widetilde{X})_q$ onto ${\cal F}$ such
that for any function $f \in{\cal F}$, (\ref{Xn+f}) and (\ref{Xn-f}) are
valid.
\end{proposition}

\smallskip

 {\bf Proof.} The uniqueness of the extension is obvious. The existence is
due to the fact that for any function $f \in{\cal F}$ and any finite subset
$M \subset q^{2{\Bbb Z}}$, there exists a polynomial $\psi \in{\Bbb C}[x]$
with $f(x)=\psi(x)$ for all $x \in M$. \hfill $\Box$

\medskip

\begin{proposition} Let $J$ be the bilateral ideal in ${\cal F}$ generated
by such elements $f \in{\Bbb F}$ that ${\tt supp}f(x)\subset q^{2{\Bbb N}}$.
Then $J$ is a submodule of the $U_q \frak{su}_{nm}$-module ${\cal F}$.
\end{proposition}

\smallskip

 {\bf Proof.} It suffices to apply (\ref{Xn+f}), (\ref{Xn-f}), and the
relations $X_j^\pm f(x)=0$ for $j \ne n$, $H_if(x)=0$ for $i=1,\ldots,N-1$.
\hfill $\Box$

\medskip

\begin{corollary} The quotient algebra ${\cal F}/J$ is a covariant
$*$-algebra.
\end{corollary}

\medskip

{\sc Remark 7.4.} The algebra of functions $f(x)$ with finite support ${\tt
supp}f \subset q^{-2{\Bbb Z}+}$ is obviously embedded into ${\cal F}/J$.
Among the above functions, a principal position is occupied by the function
$$f(x)=\left \{{1,\;x=1 \atop 0,\;x \ne 1}\right..$$
This element of ${\cal F}/J$ is denoted in the sequel by $f_0$. The next
proposition is a straightforward consequence of our definitions.

\medskip \stepcounter{theorem}

\begin{proposition}\label{tijf0}In ${\cal F}/J$, the following relations are
valid:
\begin{equation}t_{ij}f_0=0 \quad{\rm for}\quad i>m \;\&\;j>n,\end{equation}
\begin{equation}f_0t_{ij}=0 \quad{\rm for}\quad i \le m \;\&\;j \le n,
\end{equation}
\begin{equation}t_{ij}f_0=f_0t_{ij} \quad{\rm for}\quad i>m \;\&\;j \le n
\quad{\rm or}\quad i \le m \;\&\;j>n,
\end{equation}
$$xf_0=f_0x=f_0,\qquad f_0=f_0^2=f_0^*.$$
\end{proposition}

\medskip

 ($f_0$ can be treated as a q-analogue of the delta-function $\delta(x-1)$
on $\widetilde{X}$).

 One can use (\ref{Xn+f}) -- (\ref{Xn-x}) to deduce

\medskip

\begin{proposition} In the covariant $*$-algebra ${\cal F}/J$ one has
$$K_n^{\pm 1}f_0=f_0,$$
$$E_nf_0=-{q^{3/2}\over 1-q^2}t_{\{1,2,\ldots,m \}\{n,n+2,\ldots,N
\}}^{\wedge m}t^*f_0,$$
$$F_nf_0=-{q^{3/2}\over q^{-2}-1}f_0t\left(t_{\{1,2,\ldots,m
\}\{n,n+2,\ldots,N \}}^{\wedge m}\right)^*,$$ while for $j \ne n$ $$K_j^{\pm
1}f_0=f_0,\qquad E_jf_0=F_jf_0=0.$$ \end{proposition}

\medskip

 Consider a covariant $*$-algebra ${\rm Fun}(\widetilde{X})_q \subset{\cal
F}/J$ , generated by $\{t_{ij}\}_{i,j=1,\ldots,N}$ and
$f_0$.\footnote{We do not discuss in the present work whether or not this
inclusion is proper.}
Replace in the above observations the ring of polynomials ${\Bbb C}[x]$ with
the ring ${\Bbb C}[x,x^{-1}]$ to get a covariant $*$-algebra ${\rm
Fun}(\widetilde{X})_{q,x}\supset{\rm Fun}(\widetilde{X})_q$.

 The term 'finite functions' is reserved for the elements of the bilateral
ideal $D(\widetilde{X})_q$ of ${\rm Fun}(\widetilde{X})_q$ generated by
$f_0$. It is evident that
$${\rm Fun}(\widetilde{X})_q={\rm
Pol}(\widetilde{X})_q+D(\widetilde{X})_q,\qquad{\rm
Fun}(\widetilde{X})_{q,x}={\rm
Pol}(\widetilde{X})_{q,x}+D(\widetilde{X})_q,$$
and ${\rm Fun}(\widetilde{X})_q$, ${\rm Fun}(\widetilde{X})_{q,x}$,
$D(\widetilde{X})_q$ are covariant $*$-algebras.

 Define a covariant $*$-algebra ${\rm Fun}({\Bbb U})_q$ by its generators
$f_0$, $z_a^\alpha$, $a=1,\ldots,n$, $\alpha=1,\ldots,m$, and the relations
(motivated by the relations in ${\rm Fun}(\widetilde{X})_q$). The list of
relations that determine ${\rm Fun}({\Bbb U})_q$ includes all the relations
of ${\rm Pol}({\rm Mat}_{mn})_q$, and additionally
$$f_0=f_0^2=f_0^*,\qquad(z_a^\alpha)^*f_0=f_0(z_a^\alpha)=0,\quad
a=1,\ldots,n,\quad \alpha=1,\ldots,m.$$
The structure of a covariant $*$-algebra is imposed by
$$K_n^{\pm 1}f_0=f_0,\qquad E_nf_0=-{q^{1/2}\over 1-q^2}z_n^mf_0,\qquad
F_nf_0=-{q^{1/2}\over q^{-2}-1}f_0(z_n^m)^*,$$
$$(K_j^{\pm 1}-1)f_0=E_jf_0=F_jf_0=0 \quad{\rm for}\quad j \ne m.$$
(To verify the correctness of the latter definition, one has to apply
corollary \ref{Xnzaa} and the fact that $z_n^m$ quasi-commutes with
$(z_a^\alpha)^*$ for $(a,\alpha)\ne(n,m)$.) Of course, the bilateral ideal
$D({\Bbb U})_q \subset{\rm Fun}({\Bbb U})_q$ generated by $f_0$, is a
covariant $*$-algebra.

 Evidently, the map
$$i:f_0 \mapsto f_0,\qquad i:z_a^\alpha \mapsto t^{-1}t_{\{1,2,\ldots,m
\}J_{a \alpha}},\quad a=1,\ldots,n;\quad \alpha=1,\ldots,m$$
admits a unique extension up to a homomorphism of covariant $*$-algebras
$i:{\rm Fun}({\Bbb U})_q \to{\rm Fun}(\widetilde{X})_{q,x}$.

 {\sc Remark 7.7.} We have extended the structure of $U_q
\frak{sl}_N$-module algebra from ${\Bbb C}[SL_N]_q$ onto ${\rm
Fun}(\widetilde{X})_{q,x}$. In a similar way, the structure of $U_q
\frak{sl}_N^{\rm op}$-module algebra admits an extension as well (see
section 5), as one can observe from the following analogues of (\ref{Xn+f}),
(\ref{Xn-f}):
$$(H_m \otimes 1)f(x)=0,\qquad (X_m^+\otimes
1)f(x)=\frac{f(q^{-2}x)-f(x)}{q^{-2}x-x}(X_m^+\otimes 1)x,$$
$$(X_m^-\otimes
1)f(x)=(X_m^-\otimes 1)x \cdot \frac{f(q^{-2}x)-f(x)}{q^{-2}x-x},$$
$$(H_j \otimes 1)f(x)=(X_j^-\otimes 1)f(x)=0 \quad{\rm for}\quad j \ne m.$$
Hence
$$(H_m \otimes 1)f_0=0,\qquad(X_m^+\otimes 1)f_0=-{1 \over q^{-2}-1}f_0
\cdot(X_m^+\otimes 1)x,$$
$$\qquad(X_m^-\otimes 1)f_0=-{1 \over q^{-2}-1} \cdot((X_m^-\otimes
1)x)f_0.$$

 Turn to a construction of quantum analogues for homogeneous spaces $S(U_n
\times U_m)\backslash SU_{nm}$ and $(SU_n \times SU_m)\backslash SU_{nm}$.

 Remind that a vector $v \in V$ is called an $A$-invariant if $V$ is a
module over the Hopf algebra $A$ with a counit $\varepsilon$ and $a \cdot
v=\varepsilon(a)\cdot v$ for all $a \in A$.

 Let $U_q \frak{s}(\frak{gl}_m \times \frak{gl}_n)^{\rm op}\subset U_q
\frak{sl}_N^{\rm op}$ be the Hopf subalgebra generated by $K_m^{\pm 1}$,
$\{K_j^{\pm 1},E_j,F_j \}_{j \ne m}$ and $U_q \frak{sl}_m^{\rm op}\otimes
U_q \frak{sl}_n^{\rm op}\subset U_q \frak{s}(\frak{gl}_m \times
\frak{gl}_n)^{\rm op}$ -- a Hopf subalgebra generated by $\{K_j^{\pm
1},E_j,F_j \}_{j \ne m}$.

 Consider the covariant $*$-algebra ${\rm Pol}(\widetilde{X})_q$. Denote the
subalgebra of its $U_q \frak{s}(\frak{gl}_m \times \frak{gl}_n)^{\rm
op}$-invariants by ${\rm Pol}(X)_q$, and the subalgebra of $U_q
\frak{sl}_m^{\rm op}\otimes U_q \frak{sl}_n^{\rm op}$-invariants by ${\rm
Pol}(\widehat{X})_q$. Evidently, they are covariant $*$-algebras and
$${\rm Pol}(X)_q=\{f \in{\rm Pol}(\widehat{X})_q|\;(H_m \otimes 1)f=0 \}.$$
${\rm Pol}(X)_q$, ${\rm Pol}(\widehat{X})_q$ substitute the algebras of
polynomial functions on the homogeneous spaces
$$X=S(U_m \times U_n)\backslash \widetilde{X}\sim S(U_n \times
U_m)\backslash SU_{nm},$$
$$\widehat{X}=(SU_m \times SU_n)\backslash \widetilde{X}\sim(SU_n \times
SU_m)\backslash SU_{nm}.$$

 A replacement of ${\rm Pol}(\widetilde{X})_q$ in the above observations by
any of the following covariant $*$-algebras ${\rm
Pol}(\widetilde{X})_{q,x}$, ${\rm Fun}(\widetilde{X})_q$, ${\rm
Fun}(\widetilde{X})_{q,x}$, $D(\widetilde{X})_q$, allows one to produce the
covariant $*$-algebras ${\rm Pol}(X)_{q,x}\subset{\rm
Pol}(\widehat{X})_{q,x}$, ${\rm Fun}(X)_q \subset{\rm Fun}(\widehat{X})_q$,
${\rm Fun}(X)_{q,x}\subset{\rm Fun}(\widehat{X})_{q,x}$, $D(X)_q \subset
D(\widehat{X})_q$. (In every pair a smaller subalgebra is distinguished from
a larger one by the equation $(H_m \otimes 1)f=0$.)

 Note that the element $x$ is a $U_q \frak{s}(\frak{gl}_m \times
\frak{gl}_n)^{\rm op}$-invariant. Therefore the algebras ${\rm
Pol}(X)_{q,x}\subset{\rm Pol}(\widehat{X})_{q,x}$ are derivable from ${\rm
Pol}(X)_q \subset{\rm Pol}(\widehat{X})_q$ via a localization with respect
to the multiplicative system $x^{\Bbb N}$.

\bigskip

\section{Canonical isomorphism $\bf D({\Bbb U})_q \simeq D(X)_q$}

 Section 7 contains a construction of a morphism of covariant $*$-algebras
$i:{\rm Fun}({\Bbb U})_q \to{\rm Fun}(\widetilde{X})_{q,x}$. Our immediate
purpose is to prove its injectivity and to describe the image of $D({\Bbb
U})_q$ with respect to the embedding into ${\rm Fun}(\widetilde{X})_{q,x}$.

\medskip

\begin{proposition}\label{xif}The least subalgebra in ${\rm
Pol}(\widetilde{X})_{q,x}$ which contains $x$ and $if$, $f \in{\rm Pol}({\rm
Mat}_{mn})_q$, is ${\rm Pol}(X)_{q,x}$.
\end{proposition}

\smallskip

 {\bf Proof.} The least subalgebra of ${\rm Pol}(\widetilde{X})_{q,x}$ which
contains $t$, $t^{-1}$, $t^*$, $t^{*-1}$ and $if$, $f \in{\rm Pol}({\rm
Mat}_{mn})_q$, is a $U_q \frak{sl}_N$-module subalgebra. Hence it coincides
with ${\rm Pol}(\widehat{X})_{q,x}$ by a virtue of lemma \ref{ii} to be
proved below. Thus every element $\psi \in{\rm Pol}(\widehat{X})_{q,x}$ can
be written in the form $\psi=\sum \limits_{j=1}^\infty \psi_jt^j+\psi_0+\sum
\limits_{j=1}^\infty \psi_{-j}t^{*j}$ with coefficients $\{\psi_j
\}_{j=-\infty}^\infty$ from a subalgebra $F \subset{\rm Pol}(X)_{q,x}$
spanned by $x$ and $if$, $f \in{\rm Pol}({\rm Mat}_{mn})_q$. The above
representation for $\psi$ may be treated as an expansion of $\psi$ in
eigenvectors of $K_m^{\pm 1}\otimes 1$. On the other hand, ${\rm
Pol}(X)_{q,x}=\{f \in{\rm Pol}(\widehat{X})_{q,x}|\;(K_m^{\pm 1}\otimes
1)f-f=0 \}$. Thus, $\psi \in{\rm Pol}(X)_{q,x}$ implies $\psi_j=0$ for all
$j \ne 0$. That is, $\psi \in{\rm Pol}(X)_{q,x}$ implies $\psi \in F$.
\hfill $\Box$

\medskip

\begin{lemma}\label{ii}The least $U_q \frak{sl}_N$-module subalgebra of ${\rm
Pol}(\widetilde{X})_q$ which contains both $t_{\{1,2,\ldots,m
\}\{n+1,n+2,\ldots,N \}}^{\wedge m}$ and $t_{\{m+1,m+2,\ldots,N
\}\{1,2,\ldots,n \}}^{\wedge n}$, is ${\rm Pol}(\widehat{X})_{q}$.
\end{lemma}

\smallskip

 {\bf Proof.} Consider the subalgebra $\widetilde{\Bbb F}_+\subset{\Bbb
C}[SL_N]_q$ generated by $t_{ij}$ with $i \le m$. Just as in the case $q=1$,
it has simple $U_q \frak{sl}_m^{\rm op}\otimes U_q \frak{sl}_N$-isotypical
components whose generators are $t_{1N}^{a_1}\cdot \left(t_{\{1,2 \}\{N-1,N
\}}^{\wedge 2}\right)^{a_2}\cdot \ldots \cdot \left(t_{\{1,2,\ldots,m
\}\{n+1,n+2,\ldots,N \}}^{\wedge m}\right)^{a_m}$, $a_1,a_2,\ldots,a_m
\in{\Bbb Z}_+$. (A classical result of theory of invariants is applied here
(see, for example, \cite{Jo, T}), together with the fact that the
dimensionality of a simple module with highest weight $\lambda$ is
independent of $q \in(0,1]$ (see \cite{CP})). The above monomials are $U_q
\frak{sl}_m^{\rm op}$-invariant if and only if $a_1,a_2,\ldots,a_{m-1}=0$.
Hence generators of the $U_q \frak{sl}_N$-module
$$\widehat{\Bbb F}_+=\{f \in \widetilde{\Bbb F}_+|\;(E_j \otimes 1)f=(F_j
\otimes 1)f=(K_j^{\pm 1}\otimes 1)f-f=0,\;j=1,\ldots,m-1 \}$$
are the vectors $\left(t_{\{1,2,\ldots,m \}\{n+1,n+2,\ldots,N
\}}^{\wedge m}\right)^{a_m}$, $a_m \in{\Bbb Z}_+$. That is, $U_q
\frak{sl}_N$-module algebra $\widehat{\Bbb F}_+$ is generated by
$t_{\{1,2,\ldots,m \}\{n+1,n+2,\ldots,N \}}^{\wedge m}$.

 Consider the subalgebra $\widetilde{\Bbb F}_-\subset{\Bbb C}[SL_N]_q$
generated by $t_{ij}$, $i>m$, and the subalgebra $\widehat{\Bbb
F}_-=\widetilde{\Bbb F}_-\cap{\rm Pol}(\widehat{X})_q$. Just as in the case
of $\widehat{\Bbb F}_+$, one can prove that $t_{\{m+1,m+2,\ldots,N
\}\{1,2,\ldots,n \}}^{\wedge n}$ generates $\widehat{\Bbb F}_-$ as a $U_q
\frak{sl}_N$-module algebra. What remains is to apply ${\Bbb
C}[SL_N]_q=\widetilde{\Bbb F}_+\cdot \widetilde{\Bbb F}_-$, ${\rm
Pol}(\widehat{X})_q=\widehat{\Bbb F}_+\cdot \widehat{\Bbb F}_-$. \hfill
$\Box$

\medskip

\begin{theorem}\label{iiso} $i:D({\Bbb U})_q \to D(X)_q$ is an isomorphism.
\end{theorem}

\smallskip

 {\bf Proof.} Evidently, $D(\widetilde{X})_q=\widetilde{\Bbb F}_+\cdot f_0
\cdot \widetilde{\Bbb F}_-$. Thus, by proposition \ref{tijf0},
$D(\widehat{X})_q=\widehat{\Bbb F}_+\cdot f_0 \cdot \widehat{\Bbb F}_-$. On
the other hand, by a virtue of lemma \ref{ii}, $\widehat{\Bbb F}_+$ is a
subalgebra generated by the elements $t_{\{1,2,\ldots,m \}I}^{\wedge m}$,
${\tt card}(I)=m$, while $\widehat{\Bbb F}_-$ is a subalgebra generated by
the elements $t_{\{m+1,m+2,\ldots,N \}J}^{\wedge n}$, ${\tt card}(J)=n$.
Hence, in view of $xf_0=f_0x=f_0$, one has
$$D(\widehat{X})_q=\sum_{j>0}t^jD(X)_q+D(X)_q+\sum_{j>0}D(X)_qt^{*j},$$
$$D(X)_q=\left(i{\Bbb C}[{\rm Mat}_{mn}]_q \right)f_0 \left(i{\Bbb
C}[\overline{\rm Mat}_{mn}]_q \right)=iD({\Bbb U})_q.$$
That is, we have proved that $i$ is onto. To prove its injectivity,
introduce a vector space ${\cal H}={\Bbb C}[{\rm Mat}_{mn}]_q \cdot f_0
\subset D({\Bbb U})_q$, together with a representation $\Theta$ of ${\rm
Fun}({\Bbb U})_q$ in ${\cal H}$, given by $\Theta(\psi):f \mapsto \psi f$,
$\psi \in{\rm Fun}({\Bbb U})_q$, $f \in{\cal H}$.

 Equip ${\cal H}$ with such a sesquilinear form (scalar product) that
\begin{equation}\label{sp}(\psi_1f_0,\psi_2f_0)f_0=f_0
\psi_2^*\psi_1f_0,\qquad \psi_1,\psi_2 \in{\Bbb C}[{\rm Mat}_{mn}]_q.
\end{equation}
This is well defined, as one can see from $f_0 \cdot{\rm Pol}({\rm
Mat}_{mn})_q f_0={\Bbb C}f_0$.

\medskip

\begin{lemma}\label{spp}The scalar product (\ref{sp}) in ${\cal H}$ is
positive definite \footnote{That is, $(v,v)>0$ for all $v \ne 0$.} and
$\Theta$ is a $*$-representation of ${\rm Pol}({\rm Mat}_{mn})_q$ in the
pre-Hilbert space ${\cal H}$.
\end{lemma}

\smallskip

 {\bf Proof.} Remind that the space $\widetilde{\cal H}$ of the
$*$-representation $\Pi$ (see Appendix 2) is a graded vector space:
$\widetilde{\cal H}=\bigoplus \limits_{j=0}^\infty \widetilde{\cal H}_j$,
and $\widetilde{\cal H}_0={\Bbb C}\cdot {\bf e}_0$. For all $v \in
\widetilde{\cal H}$ and all $a,\alpha$, one has
$${\rm deg}(\Pi(z_a^\alpha)v)={\rm deg}\,v+1,\qquad{\rm
deg}(\Pi(z_a^\alpha)^*v)={\rm deg}\,v-1$$
by a virtue of (\ref{tijfx}). Hence, $(\Pi(\psi){\bf e}_0,{\bf e}_0)f_0=f_0
\psi f_0$ for all $\psi \in{\rm Pol}({\rm Mat}_{mn})_q$. Thus, the map
$$j:{\cal H}\to \widetilde{\cal H},\qquad j:\Theta(\psi)f_0 \mapsto
\Pi(\psi){\bf e}_0,\qquad \psi \in{\Bbb C}[{\rm Mat}_{mn}]_q$$
is well defined and intertwines the scalar products.

 Now apply propositions A2.2.3 and \ref{xif} to conclude that the map $j$ is
onto. Hence $j({\Bbb C}[{\rm Mat}_{mn}]_{qk}\cdot f_0)=\widetilde{\cal H}_k$
for all $k \in{\Bbb Z}_+$. On the other hand, ${\rm dim}\,\widetilde{\cal
H}_k={\rm dim}\,{\Bbb C}[{\rm Mat}_{mn}]_{qk}$. Therefore, $j$ is
one-to-one, and the representations $\Theta$ and $\Pi$ are unitarily
equivalent. \hfill $\Box$

\medskip \stepcounter{theorem}

{\sc Remark 8.5.} There exists a unique extension of $\widetilde{\Pi}$ onto
${\rm Fun}(\widetilde{X})_q$ such that $\widetilde{\Pi}(f_0)$ is the
projection onto $\widetilde{\cal H}_0$ with kernel $\bigoplus \limits_{k \ne
0}\widetilde{\cal H}_k$. One can observe from the proof of lemma \ref{spp}
that the representations $\Theta$ and $\Pi=\widetilde{\Pi}\circ i$ of ${\rm
Fun}({\Bbb U})_q$ are unitarily equivalent.

\medskip

 Turn back to proving the injectivity of $i$. By a virtue of remark 8.5, it
suffices to prove that the homomorphism $\Theta:D({\Bbb U})_q \to{\rm
End}({\cal H})$ is an embedding. The linear map
$$m:{\Bbb C}[{\rm Mat}_{mn}]_q \cdot f_0 \otimes f_0 \cdot{\Bbb
C}[\overline{\rm Mat}_{mn}]_q \to D({\Bbb U})_q,\qquad m:f_1 \otimes f_2
\mapsto f_1f_2$$
is one-to-one, as one can easily deduce via an application of the well known
diamond lemma \cite{B} to producing monomial bases in the vector spaces
${\Bbb C}[{\rm Mat}_{mn}]_q \cdot f_0$, $f_0 \cdot{\Bbb C}[\overline{\rm
Mat}_{mn}]_q$, $D({\Bbb U})_q$. Thus,
$$D({\Bbb U})_q \simeq({\Bbb C}[{\rm Mat}_{mn}]_q \cdot
f_0)\otimes(f_0 \cdot{\Bbb C}[\overline{\rm Mat}_{mn}]_q) \simeq{\cal
H}\otimes{\cal H}^*,$$
and the representation $\Theta:D({\Bbb U})_q \to{\rm End}({\cal H})$ reduces
to the canonical linear map ${\cal H}\otimes{\cal H}^*\to{\rm End}({\cal
H})$. What remains is to observe that this map is an embedding. The theorem
\ref{iiso} is proved. \hfill $\Box$

\medskip \stepcounter{theorem}

{\sc Remark 8.6.} One can easily deduce a description of the image
$i(D({\Bbb U})_q)$ in ${\rm End}({\cal H})$. Equip ${\cal H}$ with a
gradation
$${\cal H}=\bigoplus_{k=0}^\infty{\cal H}_k,\qquad{\cal H}_k={\Bbb C}[{\rm
Mat}_{mn}]_{q,k}\cdot f_0,$$
and let ${\rm End}({\cal H})_f$ stand for the algebra of finite dimensional
finite degree operators in ${\cal H}$. Since ${\tt dim}\,{\cal H}_k<\infty$,
$k<\infty$, one has ${\rm End}({\cal H})_f \simeq{\cal H}\otimes{\cal H}^*$.
Hence, $\Theta$ provides a 'canonical' isomorphism of algebras $D({\Bbb
U})_q \to{\rm End}({\cal H})_f$.

\medskip

 Let $P_k$ be the projection  in ${\cal H}$ onto the homogeneous component
${\cal H}_k={\Bbb C}[{\rm Mat}_{mn}]_{q,k}\cdot f_0$ with kernel $\bigoplus
\limits_{j \ne k}{\Bbb C}[{\rm Mat}_{mn}]_{q,j}\cdot f_0$, and let ${\Bbb
C}[{\rm Mat}_{mn}]_{q,k}\cdot{\Bbb C}[\overline{\rm Mat}_{mn}]_{q,-l}$ be
the linear span of
$$\{f_+\cdot f_-\in{\rm Pol}({\rm Mat}_{mn})_q|\;f_+\in{\Bbb C}[{\rm
Mat}_{mn}]_q,\;f_-\in{\Bbb C}[\overline{\rm Mat}_{mn}]_q \}.$$

\medskip

\begin{lemma}\label{oo}For all $k,l \in{\Bbb Z}_+$, the map
$${\Bbb C}[{\rm Mat}_{mn}]_{q,k}\cdot{\Bbb C}[\overline{\rm
Mat}_{mn}]_{q,-l}\to{\rm Hom}({\cal H}_l,{\cal H}_k);\qquad f \mapsto P_k
\Theta(f)|_{{\cal H}_l}$$
is one-to-one.
\end{lemma}

\smallskip

 {\bf Proof.} Both ${\cal H}_k$, ${\cal H}_l$ are finite dimensional Hilbert
spaces. Arguing just as in the proof of theorem \ref{iiso}, we get
$${\Bbb C}[{\rm Mat}_{mn}]_{q,k}\cdot{\Bbb C}[\overline{\rm
Mat}_{mn}]_{q,-l}\simeq{\Bbb C}[{\rm Mat}_{mn}]_{q,k}\otimes{\Bbb
C}[\overline{\rm Mat}_{mn}]_{q,-l}\simeq$$
$$\simeq{\Bbb C}[{\rm Mat}_{mn}]_{q,k}\cdot f_0 \otimes f_0 \cdot{\Bbb
C}[\overline{\rm Mat}_{mn}]_{q,-l}\simeq{\cal H}_k \otimes{\cal H}_l^*
\simeq{\rm Hom}({\cal H}_l,{\cal H}_k).$$
What remains is to use the fact that the resulting linear map
$${\Bbb C}[{\rm Mat}_{mn}]_{q,k}\cdot{\Bbb C}[\overline{\rm
Mat}_{mn}]_{q,-l}\to_{_{_{\!\!\!\!\!\!\!\!\textstyle \sim}}}{\rm
Hom}({\cal H}_l,{\cal H}_k)$$
coincides with the operator described in the statement of this lemma. (In
fact, let $f=f_+f_-^*$, $f_+\in{\Bbb C}[{\rm Mat}_{mn}]_{q,k}$, $f_-\in{\Bbb
C}[{\rm Mat}_{mn}]_{q,l}$. Then for all $\psi_+\in{\Bbb C}[{\rm
Mat}_{mn}]_{q,k}$, $\psi_-\in{\Bbb C}[{\rm Mat}_{mn}]_{q,l}$ one has
$$(\psi_+f_0,P_k
\Theta(f_+f_-^*)\psi_-f_0)=(\psi_+f_0,f_+f_-^*\psi_-f_0)=
(f_+^*\psi_+f_0,f_-^*\psi_-f_0)=$$
$$=((\psi_+f_0,f_+f_0)f_0,(\psi_-f_0,f_-f_0)f_0)=(\psi_+f_0,f_+f_0)\cdot
\overline{(\psi_-f_0,f_-f_0)}.\eqno \Box$$

\medskip

\begin{proposition}\label{em}The homomorphism $\Theta:{\rm Pol}({\rm
Mat}_{mn})_q \to{\rm End}({\cal H})$ is an embedding.
\end{proposition}

\smallskip

 {\bf Proof.} Equip the vector space ${\rm Pol}({\rm Mat}_{mn})_q$ with a
bigradation
$${\rm Pol}({\rm Mat}_{mn})_q=\bigoplus_{k,l=0}^\infty{\Bbb C}[{\rm
Mat}_{mn}]_{q,k}\cdot{\Bbb C}[\overline{\rm Mat}_{mn}]_{q,-l}$$
(as one can easily verify, this is well defined). We need also a standard
partial order relation on ${\Bbb Z}_+^2$:
$$(k_1,l_1)\le(k_2,l_2)\qquad \Leftrightarrow \qquad k_1 \le k_2 \quad \&
\quad l_1 \le l_2.$$

 Assume that our statement is wrong and $\Theta(f)=0$ for some $f \in{\rm
Pol}({\rm Mat}_{mn})_q$, $f \ne 0$. Consider a homogeneous component
$f_{kl}$ of $f$ with minimal bidegree $(k,l)$. (Such homogeneous component
certainly exists, but it is not unique for a given $f \in{\rm Pol}({\rm
Mat}_{mn})_q$). Let ${\cal H}_j={\Bbb C}[{\rm Mat}_{mn}]_{q,j}\cdot f_0$,
and $P_k:{\cal H}\to{\cal H}_k$ be the projection onto ${\cal H}_k$ with
kernel $\bigoplus \limits_{j \ne k}{\cal H}_j$. Since $f_{kl}$ is of a
minimal bidegree, one has $P_k \Theta(f_{kl})|_{{\cal H}_l}=P_k
\Theta(f)|_{{\cal H}_l}=0$, $f_{kl}\ne 0$, which contradicts the statement
of lemma \ref{oo}. \hfill $\Box$

\medskip

\begin{corollary} i) The morphism $i:{\rm Pol}({\rm Mat}_{mn})_q \to{\rm
Fun}(\widetilde{X})_{q,x}$ is an embedding.\\
ii) The restriction of $\widetilde{\Pi}$ onto the subalgebra ${\rm
Pol}(X)_{q,x}$ is a faithful $*$-representation of this subalgebra.
\end{corollary}

\smallskip

 {\bf Proof.} The first statement follows from the equivalence of $\Theta$
and $\Pi=\widetilde{\Pi}\circ i$, and the faithfulness of $\Theta$
established in proposition \ref{em}.

 Now turn to proving the second statement. Suppose that $\psi \in{\rm
Pol}(X)_{q,x}$ and $\widetilde{\Pi}(\psi)=0$. By proposition \ref{xif} and
the relations (\ref{tijfx}), there exist such elements
$\psi_1,\psi_2,\ldots,\psi_M \in{\rm Pol}({\rm Mat}_{mn})_q$ that $\psi=\sum
\limits_{k=0}^Mi(\psi_k)x^k$. In \cite{SSV7} an element $y \in{\rm Pol}({\rm
Mat}_{mn})_q$ is found with the property $iy=x^{-1}$. Hence
$\psi=i(\Psi)x^M$ with $\Psi=\sum \limits_{k=0}^M \psi_ky^{M-k}$. It follows
from $\widetilde{\Pi}(\psi)=0$ that $\widetilde{\Pi}i(\Psi)\cdot
\left(\widetilde{\Pi}(x)\right)^M=0$. Observe that $\widetilde{\Pi}(x)$ is
invertible, and so $\widetilde{\Pi}i(\Psi)=0$, $\Psi=0$, $\psi=0$. \hfill
$\Box$

\medskip

\begin{proposition}The morphism of covariant algebras $i:{\rm Fun}({\Bbb
U})_q \to{\rm Fun}(\widetilde{X})_{q,x}$ is an embedding.
\end{proposition}

\smallskip

 {\bf Proof.} It was proved before that $i$ is an embedding while restricted
onto the subalgebras $D({\Bbb U})_q$ and ${\rm Pol}({\rm Mat}_{mn})_q$.

 Let $i(f_1+f_2)=0$, $f_1 \in D({\Bbb U})_q$, $f_2 \in{\rm Pol}({\rm
Mat}_{mn})_q$. By a virtue of Remark 8.6, $\Theta(f_1){\cal H}\subset
\bigoplus \limits_{j=0}^{M-1}{\cal H}_j$ for some $M \in{\Bbb N}$. It
follows that $\Theta(f_2){\cal H}\subset \bigoplus \limits_{j=0}^{M-1}{\cal
H}_j$. Hence all the elements of ${\Bbb C}[\overline{\rm
Mat}_{mn}]_{q,-M}f_2$ are in the kernel of $\Theta$. By proposition
\ref{em}, ${\Bbb C}[\overline{\rm Mat}_{mn}]_{q,-M}f_2=0$. We claim this
implies $f_2=0$. In fact, the invertibility of the linear maps
$R_{\overline{U}U}$, $R_{\overline{V}V}$ for all $q \in(0,1)$ allows one to
apply diamond lemma to prove that $${\rm Pol}({\rm
Mat}_{mn})_q=\bigoplus_{k,l=0}^\infty{\Bbb C}[\overline{\rm
Mat}_{mn}]_{q,-l}\cdot{\Bbb C}[{\rm Mat}_{mn}]_{q,k},$$ via producing bases
of lexicographically ordered monomials in each of the subspaces ${\Bbb
C}[\overline{\rm Mat}_{mn}]_{q,-l}$, ${\Bbb C}[{\rm Mat}_{mn}]_{q,k}$. If
$\psi$ is the first (lowest) element of such basis in ${\Bbb
C}[\overline{\rm Mat}_{mn}]_{q,-M}$, then obviously $\psi f_2=0$ implies
$f_2=0$. \hfill $\Box$

\bigskip

\section{An invariant integral}

 Consider the Hopf subalgebra $U_q \frak{p}_+\subset U_q \frak{sl}_N$
generated by $K_n^{\pm 1}$, $E_n$, and $K_j^{\pm 1}$, $E_j$, $F_j$, $j \ne
n$. By a virtue of the relation $E_nf_0=-{\textstyle q^{1/2}\over \textstyle
1-q^2}z_n^mf_0$ from section 7, one has

\medskip

\begin{proposition} ${\cal H}$ is a $U_q \frak{p}_+$-submodule of the $U_q
\frak{p}_+$-module $D({\Bbb U})_q$.
\end{proposition}

\medskip

 Let $\Gamma$ stand for the associated representation of $U_q \frak{p}_+$ in
${\cal H}$.

 We use in what follows the standard scalar product in the Cartan subalgebra
$\frak{h}$:
$$(H_i,H_j)=\left \{\begin{array}{ccl}2 &,& i=j \\ -1 &,& |i-j|=1 \\ 0 &,&
{\rm otherwise}\end{array}\right.$$
and the element $\check{\rho}\in \frak{h}$ given by $(\check{\rho},H_i)=1$,
$i=1,\ldots,N-1$.

\medskip

\begin{theorem}\label{iin} The linear functional
$$\nu:D({\Bbb U})_q \to{\Bbb C},\qquad \int \limits_{{\Bbb U}_q}fd
\nu \stackrel{\rm def}{=}{\tt tr}(\Theta(f)\Gamma(e^{h \check{\rho}}))$$
is well defined, $U_q \frak{sl}_N$-invariant and positive
\footnote{Positive in the sense that $\int \limits_{{\Bbb U}_q}fd \nu>0$
for all non-zero non-negative elements of the $*$-algebra
$D({\scriptstyle \Bbb U})_q$.}.
\end{theorem}

\smallskip

 {\bf Proof.} One can deduce that $\nu$ is well defined and positive from
the results of section 8 since ${\cal H}$ is a pre-Hilbert space, the
$*$-representation $\Theta$ is faithful, and $\Theta(f)$, $f \in D({\Bbb
U})_q$, are finite dimensional finite degree operators. The proof of $U_q
\frak{sl}_N$-invariance of $\nu$ is just the same as that in the special
case $m=n=1$ \cite{SSV2}.

 Specifically, by a virtue of (\ref{*H}) for $F=D({\Bbb U})_q$, $U_q
\frak{sl}_N$-invariance of this integral follows from its $U_q
\frak{p}_+$-invariance and its realness: $\int \limits_{{\Bbb U}_q}f^*d
\nu=\overline{\int \limits_{{\Bbb U}_q}fd \nu}$, $f \in D({\Bbb U})_q$. So
what remains is to prove the $U_q \frak{p}_+$-invariance of the integral we
deal with. It follows from the covariance of $D({\Bbb U})_q$ that the linear
map $D({\Bbb U})_q \otimes{\cal H}\to{\cal H}$, $f \otimes v \mapsto fv$, $f
\in D({\Bbb U})_q$, $v \in{\cal H}$, is a morphism of $U_q
\frak{p}_+$-modules. Hence the associated linear map $D({\Bbb U})_q \to{\cal
H}\otimes{\cal H}^*$ is also a morphism of $U_q \frak{p}_+$-modules (see
\cite[proposition 1.2]{SSV2}). So one needs to use the fact that the square
of the antipode $S$ is an inner automorphism $S^2(\xi)=e^{h
\check{\rho}}\cdot \xi \cdot e^{-h \check{\rho}}$, $\xi \in U_q
\frak{sl}_N$, and to apply the general argument given below (it is well
known from the theory of Hopf algebras \cite{CP}) to the $U_q
\frak{p}_+$-module ${\cal H}$.

 Let $A$ be a Hopf algebra and $\Gamma$ its representation in a vector space
$V$. Then $V$, $V^*$, $V^{**},\ldots$ are $A$-modules, while the standard
embedding $i_0:V \hookrightarrow V^{**}$, is not in general a morphism of
$A$-modules (unless $S^2={\tt id}$). Let $u \in A$ be such that $S^2(\xi)=u
\cdot \xi \cdot u^{-1}$ for all $\xi \in A$. Then the embedding $i_1=i_0
\Gamma(u):V \hookrightarrow V^{**}$ is a morphism of $A$-modules since
$i_0S^2(\xi)v=\xi i_0 v$ for all $v \in V$, $\xi \in A$.

 We observe that the composition of the linear map $i_1 \otimes{\tt id}:V
\otimes V^*\to V^{**}\otimes V^*$ and the canonical pairing $V^{**}\otimes
V^*\to{\Bbb C}$ is a morphism of $A$-modules, i.e. an invariant integral
\index{invariant integral}. This invariant integral can be written in the
form ${\tt tr}_q(A)={\tt tr}(A \Gamma(u))$, $A \in V \otimes V^*\subset{\rm
End}_{\Bbb C}(V)$ via an application of the canonical embedding $V \otimes
V^*\hookrightarrow{\rm End}_{\Bbb C}(V)$. \hfill $\Box$

\medskip

 To conclude, we apply theorem \ref{iin} for producing a positive invariant
integral on the quantum principal homogeneous space.

 A passage from functions on $\widetilde{X}$ to functions on $X$ could be
done via averaging with respect to an action of the compact group $S(U_m
\times U_n)$. We do this in the quantum case.

 Consider the bilateral ideal $J \subset{\Bbb C}[SL_N]_q$ generated by
$t_{kl}$ with $k \le m \quad \& \quad l>m$ or $k>m \quad \& \quad l \le
m$, and the canonical onto morphism $$j:{\Bbb C}[SL_N]_q \to{\Bbb C}[S(GL_m
\times GL_n)]_q,$$ with ${\Bbb C}[S(GL_m \times GL_n)]_q={\Bbb
C}[SL_N]_q/J$.  Introduce the notation $${\Bbb C}[S(U_m \times
U_n)]_q=({\Bbb C}[S(GL_m \times GL_n)]_q,\star)$$ for the 'algebra of
regular functions on the compact quantum group $S(U_m \times U_n)_q$'
\footnote{It is easy to prove that $J^\star \subset J$. For example,
obviously, $t_{1N}^\star \in J$, $t_{N1}^\star \in J$, and for other
generators of $J$ the covariance argument is applicable.}.

\medskip

\begin{lemma} The composition $\widetilde{\Delta}$ of homomorphisms
$\Delta:{\Bbb C}[SL_N]_q \to{\Bbb C}[SL_N]_q \otimes{\Bbb C}[SL_N]_q$, $j
\otimes{\tt id}:{\Bbb C}[SL_N]_q \otimes{\Bbb C}[SL_N]_q \to{\Bbb C}[S(GL_m
\times GL_n)]_q \otimes{\Bbb C}[SL_N]_q$ is a homomorphism of $*$-algebras
$\widetilde{\Delta}:{\rm Pol}(\widetilde{X})_q \to{\Bbb C}[S(U_m \times
U_n)]_q \otimes{\rm Pol}(\widetilde{X})_q$.
\end{lemma}

\smallskip

 {\bf Proof.} By the definition of involution in the $*$-algebra ${\Bbb
C}[SU_N]_q$, $j \otimes{\tt id}(\Delta)$ is a homomorphism of $*$-algebras
${\Bbb C}[SU_N]_q \to{\Bbb C}[S(U_m \times U_n)]_q \otimes{\Bbb C}[SU_N]_q$.
What remains is to apply the relations (\ref{*Tij}), (\ref{tijst}). \hfill
$\Box$

\medskip

 Extend $\widetilde{\Delta}$ up to a homomorphism of $*$-algebras
$\widetilde{\Delta}:{\rm Fun}(\widetilde{X})_q \to{\Bbb C}[S(U_m \times
U_n)]_q \otimes{\rm Fun}(\widetilde{X})_q$ via $\widetilde{\Delta}f_0=1
\otimes f_0$. (The existence and uniqueness of such extension follows from
the definitions of $f_0$ and the $*$-algebra ${\rm Fun}(\widetilde{X})_q$).

 Let $\widetilde{\mu}:{\Bbb C}[S(U_m \times U_n)]_q \to{\Bbb C}$ be the
invariant integral on the 'compact quantum group $S(U_m \times U_n)_q$'
normalized by $\displaystyle \int \limits_{S(U_m \times U_n)_q}1d
\widetilde{\mu}=1$ \cite{CP}, and $\widetilde{\nu}:D(X)_q \to{\Bbb C}$ an
invariant integral transferred from $\nu:D({\Bbb U})_q \to{\Bbb C}$,
$\displaystyle \int \limits_{{\Bbb U}_q}fd \nu={\tt tr}(\Theta(f)\Gamma(e^{h
\check{\rho}}))$ via the canonical isomorphism $D({\Bbb U})_q \simeq
D(X)_q$.

\medskip

\begin{proposition} The linear functional $(\widetilde{\mu}\otimes
\widetilde{\nu})\widetilde{\Delta}:D(\widetilde{X})_q \to{\Bbb C}$ is
positive and $U_q \frak{sl}_N$-invariant.
\end{proposition}

\smallskip

 {\bf Proof.} The scalar product $(f_1,f_2)=\widetilde{\mu}\otimes
\widetilde{\nu}(f_2^*f_1)$ in ${\Bbb C}[S(U_m \times U_n)]_q \otimes D(X)_q$
is positive definite, as one can see from theorem \ref{iin} and the
orthogonality relations for a compact quantum group \cite{CP}. Hence
$\widetilde{\mu}\otimes \widetilde{\nu}(f^*f)>0$ for $f \ne 0$, and the
positivity of the linear functional $(\widetilde{\mu}\otimes
\widetilde{\nu})\widetilde{\Delta}$ now follows from the injectivity of
$\widetilde{\Delta}:D(\widetilde{X})_q \to{\Bbb C}[S(U_m \times U_n)]_q
\otimes D(\widetilde{X})_q$. The $U_q \frak{sl}_N$-invariance follows from
the fact that the 'averaging operator' $(\widetilde{\mu}\otimes
id)\widetilde{\Delta}D(\widetilde{X})_q \to D(X)_q$ is a
morphism of $U_q \frak{sl}_N$-modules. \hfill $\Box$

\bigskip

\section{Concluding notes}

 A number of commutation relations obtained in sections 1 -- 6 are known
within a different approach to function theory on quantum complex manifolds
\cite{CHZ}. The conceptual equivalence of both approaches is due to the
isomorphism of quantum homogeneous spaces ${\Bbb U}$ and $X$ (see section
8).

 There is a wide class of quantum homogeneous spaces for which our explicit
formula for invariant integral is plausible. All the related covariant
algebras are derivable via a factorization from the quantum universal
enveloping algebra $U_q \frak{sl}_N$ equipped with the adjoint action
\cite{CP}. As one can observe already in the case of quantum disc (see
\cite{SSV5}), the algebras of functions considered in the present work
\underline{are not in the above wide class of covariant algebras}, although
being derivable from those by passage to a limit. That kind of passage to a
limit was investigated  before by Berezin within his approach to
quantization of bounded symmetric domains \cite{Ber1, Ber2}.

 In the first six sections of the present work, ${\Bbb C}(q^{1/s})$ worked
as a ground field, with $s$ being a natural number whose value was not
specified precisely. It follows from the subsequent descriptions of the
covariant algebras ${\rm Pol}({\rm Mat}_{mn})_q$, ${\rm
Fun}(\widetilde{X})_q$ (in terms of their generators and relations) that the
ground field could be chosen to be ${\Bbb C}(q^{1/2})$.

\bigskip

\section*{Appendix 1. Universal R-matrix and quantum Weyl group}

 The subject of this appendix is to remind some well known results of
quantum group theory. We follow S. Levendorskii and Ya. Soibelman \cite{LS1,
LS2}. A large part of these results were independently obtained by A.
Kirillov and N. Reshetikhin \cite{KR, CP}. A more general and rather
complete exposition of the background on quantum group theory can be
found in a remarkable surway of M. Rosso \cite{R}.

 Consider a reduced decomposition $w_0=s_{i_1}\cdot s_{i_2}\ldots s_{i_M}$,
$M=N(N-1)/2$, of the longest permutation $w_0=(N,N-1,\ldots,2,1)\in S_N$.
Our purpose is to associate to each such reduced decomposition a linear
order relation on the set of positive roots of the Lie algebra
$\frak{sl}_N$, and then a basis in the vector space $U_q \frak{sl}_N$.
Remind also that the simple roots $\alpha_i$, $i=1,\ldots,N-1$, are
given by $\alpha_i(H_i)=a_{ij}$, $i,j=1,\ldots,N-1$, with $(a_{ij})$ being
the Cartan matrix (\ref{Cm}). We use the following linear order relation on
the set of positive roots:
$$\beta_1=\alpha_1,\qquad \beta_2=s_{i_1}(\alpha_{i_2}),\qquad
\beta_3=s_{i_1}s_{i_2}(\alpha_{i_3}),\qquad \ldots
\qquad,\beta_M=s_{i_1}\ldots s_{i_{M-1}}(\alpha_{i_M}).$$

 The work \cite{LS1} associates to the generators $S_i$, $i=1,\ldots,N-1$,
of the Weyl group the automorphisms $T_i$ of $U_q \frak{sl}_N$, which
differ inessentially from Lusztig automorphisms (see \cite{CP},\cite{R}).
Note that, in particular, $$T_i(K_j)=\left \{\begin{array}{ccl}K_j^{-1} &,&
i=j \\ K_iK_j &,& |i-j|=1 \\ K_j &,& {\rm
otherwise}\end{array}\right..\eqno(A1.1)$$

  Associate to each simple root $\alpha_i$ the
generators $E_i$, $F_i$ of $U_q\frak{sl}_N$. The above map defined on the
 family of simple roots is extendable onto the set of all positive roots:
 $E_{\beta_s}=T_{i_1}T_{i_2}\ldots T_{i_{s-1}}(E_{i_s})$,
$F_{\beta_s}=T_{i_1}T_{i_2}\ldots T_{i_{s-1}}(F_{i_s})$.

\medskip

\noindent{\bf Proposition A1.1.} {\it $E_{\beta_1}^{k_1}\cdot
E_{\beta_2}^{k_2}\cdot \ldots \cdot E_{\beta_M}^{k_M}$,
$(k_1,k_2,\ldots,k_M)\in{\Bbb Z}_+^M$, constitute a basis in the vector
space $U_q \frak{N}_+$; $F_{\beta_1}^{j_1}\cdot F_{\beta_2}^{j_2}\cdot
\ldots \cdot F_{\beta_M}^{j_M}$, $(j_1,j_2,\ldots,j_M)\in{\Bbb Z}_+^M$,
constitute a basis in the vector space $U_q \frak{N}_-$;
$F_{\beta_1}^{k_1}\cdot F_{\beta_2}^{k_2}\cdot \ldots F_{\beta_M}^{k_M}\cdot
K_1^{i_1}\cdot K_2^{i_2}\cdot \ldots \cdot K_{N-1}^{i_{N-1}}\cdot
E_{\beta_1}^{j_1}\cdot E_{\beta_2}^{j_2}\cdot \ldots \cdot
E_{\beta_M}^{j_M}$, $(k_1,k_2,\ldots,k_M),\;(j_1,j_2,\ldots,j_M)\in{\Bbb
Z}_+^M$, $(i_1,i_2,\ldots,i_{N-1})\in{\Bbb Z}^{N-1}$, constitute a basis in
the vector space $U_q \frak{sl}_N$.}

\medskip

 Let $\Theta$ be the antiautomorphism of $U_q \frak{sl}_N$, given by
$\Theta(E_i)=F_i$, $\Theta(F_i)=E_i$, $\Theta(K_i^{\pm 1})=K_i^{\mp 1}$,
$i=1,\ldots,N-1$, and let also
$\widetilde{F}_{\beta_i}=\Theta(E_{\beta_i})$.

\medskip

\noindent{\bf Corollary A1.2.} {\it $\widetilde{F}_{\beta_M}^{k_M}\cdot
\widetilde{F}_{\beta_{M-1}}^{k_{M-1}}\cdot \ldots \cdot F_{\beta_1}^{k_1}$,
$(k_1,\ldots,k_M)\in{\Bbb Z}_+^M$, constitute a basis in the vector
space $U_q \frak{N}_-$, and $\widetilde{F}_{\beta_M}^{k_M}\cdot
\widetilde{F}_{\beta_{M-1}}^{k_{M-1}}\cdot \ldots
\widetilde{F}_{\beta_1}^{k_1}\cdot E_{\beta_1}^{j_1}\cdot
E_{\beta_2}^{j_2}\cdot \ldots \cdot E_{\beta_M}^{j_M}\cdot
K_1^{i_1}\cdot K_2^{i_2}\cdot \ldots \cdot
K_{N-1}^{i_{N-1}}$, $(k_1,k_2,\ldots,k_M),\;(j_1,j_2,\ldots,j_M)\in{\Bbb
Z}_+^M$, $(i_1,i_2,\ldots,i_{N-1})\in{\Bbb Z}^{N-1}$, constitute a basis in
the vector space $U_q \frak{sl}_N$.}

\medskip

 We are about to apply corollary A1.2 to constructing the bases of the
vector spaces $V_-(\lambda)$. For that, we use the notation of section 2,
together with the class of reduced decompositions for the element $w_0$
described there.

 Equip $U_q \frak{sl}_N$ with the gradation (cf. section 2): ${\rm
deg}(E_n)=1$, ${\rm deg}(F_n)=-1$, ${\rm deg}(K_n^{\pm 1})=0$;\  ${\rm
deg}(E_j)={\rm deg}(F_j)={\rm deg}(K_j^{\pm 1})=0$ for $j \ne n$. The
automorphisms $T_i$, $i \ne n$, preserve this gradation since the latter is
determined by the element $H_0$ defined in section 2. Hence
$$T_i(U_q \frak{sl}_n \otimes U_q \frak{sl}_m)=U_q \frak{sl}_n \otimes U_q
\frak{sl}_m,\qquad i \ne n,$$
since $T_iU_q \frak{N}_\pm \subset U_q \frak{N}_\pm$, and
$$U_q \frak{sl}_n \otimes U_q \frak{sl}_m \cap U_q \frak{N}_\pm=\{\xi \in
U_q \frak{N}_\pm|\;{\rm deg}(\xi)=0 \}.\eqno(A1.3)$$

 Impose the notation $M'=M-mn={\textstyle m(m-1)\over \textstyle
2}+{\textstyle n(n-1)\over \textstyle 2}$. It follows from (A1.3) that ${\rm
deg}(\widetilde{F}_{\beta_j})=0$ for $j \le M'$. Just in the same way as in
the case $q=1$ one deduces that ${\rm deg}(\widetilde{F}_{\beta_j})\in
\{-1,0 \}$. Thus ${\rm deg}(\widetilde{F}_{\beta_j})=-1$ for $j>M'$, and
$${\rm deg}(\widetilde{F}_{\beta_M}^{k_M}\ldots
\widetilde{F}_{\beta_1}^{k_1})=-\sum_{j=M'+1}^{M}k_j.\eqno(A1.4)$$
Now it follows from (A1.4) that the elements
$$\widetilde{F}_{\beta_{M'}}^{k_{M'}}\widetilde{F}_{\beta_{M'-1}}^{k_{M'-1}}
\ldots \widetilde{F}_{\beta_1}^{k_1},\qquad (k_1,\ldots,k_{M'})\in{\Bbb
Z}_+^{M'},$$
constitute a basis in the vector space $U_-=U_q \frak{N}_- \cap(U_q
\frak{sl}_n \otimes U_q \frak{sl}_m)$, while
$$\widetilde{F}_{\beta_M}^{k_M}\widetilde{F}_{\beta_{M-1}}^{k_{M-1}} \ldots
\widetilde{F}_{\beta_1}^{k_1},\qquad {\rm with}\qquad
\sum_{j=1}^{M'}k_j>0, \eqno(A1.5)$$ form a basis in $U_q \frak{N}_- \cdot
U_-$. Therefore, the vectors (A1.5) form a basis in the kernel of the linear
map $U_q \frak{N}_-\to V_-(0)$, $\xi \mapsto \xi v_-(0)$. Thus the vectors
$$\widetilde{F}_{\beta_M}^{k_M}\widetilde{F}_{\beta_{M-1}}^{k_{M-1}} \ldots
\widetilde{F}_{\beta_{M'+1}}^{k_{M'+1}}\cdot v_-(0),\qquad
(k_{M'+1},\ldots,k_M)\in{\Bbb Z}_+^{mn},$$
constitute a basis in the vector space $V_-(0)$. By a virtue of (A1.4),
$${\rm
deg}(\widetilde{F}_{\beta_M}^{k_M}\widetilde{F}_{\beta_{M-1}}^{k_{M-1}}
\ldots
\widetilde{F}_{\beta_{M'+1}}^{k_{M'+1}}\cdot v_-(0))=-\sum_{j=M'+1}^{M}k_j.$$

\medskip

 {\sc Remark A1.3.} In section 2 the Hopf algebra $U_A \subset U_q
\frak{sl}_N$ was considered over the ring $A={\Bbb C}[q^{1/s},q^{-1/s}]$. It
follows from the definition of the automorphisms $T_i$ (see \cite{LS1}) that
$T_iU_A=U_A$ for all $i \ne n$. Hence all the basis vectors (A1.2) of the
vector space $U_q \frak{sl}_N$ are in the lattice $U_A$.

\medskip

 Let $V_1$, $V_2$ be $U_q \frak{sl}_N$-modules satisfying the integrity
condition for weights as in section 2. Our additional assumption is that
either $V_1$ possesses a highest weight or $V_2$ possesses a lowest weight.
Under a suitable choice of the ground field ${\Bbb C}(q^{1/s})$ the formula
below determines a linear operator $R_{V_1,V_2}$ in $V_1 \otimes V_2$:
$$R=\exp_{q^2}\left((q^{-1}-q)E_{\beta_1}\otimes F_{\beta_1}\right)
\cdot \ldots \cdot
\exp_{q^2}\left((q^{-1}-q)E_{\beta_M}\otimes
F_{\beta_M}\right)q^{-t_0},\eqno(A1.6)$$
with $\exp_{q^2}(u)=\displaystyle \sum \limits_{k=0}^\infty{u^k \over
(k)_{q^2}!}$;
$$(k)_{q^2}!=\prod \limits_{j=1}^k{1-q^{2j}\over 1-q^2},\eqno(A1.7)$$
$t_0=\sum \limits_{i,j=1}^{N-1}c_{ij}H_i \otimes H_j$, and
$(c_{ij})_{i,j=1,\ldots,N-1}$ being the inverse matrix with respect to the
Cartan matrix $(a_{ij})_{i,j=1,\ldots,N-1}$.

 Now we use the relation $\alpha_i(H_j)=a_{ij}$, $i,j=1,\ldots,N-1$, to get
a different description of $t_0$:
$$\alpha_i \otimes \alpha_j(t_0)=a_{ij},\qquad i,j=1,\ldots,N-1.$$
Also, an application of the bilinear scalar product $(H_i,H_j)=a_{ij}$,
$i,j=1,\ldots,N-1$, in the Cartan subalgebra, we get the third description
of $t_0$:
$$(t_0,H_i \otimes H_j)=(H_i,H_j);\qquad i,j=1,\ldots,N-1.$$
That is, $t_0=\displaystyle \sum \limits_{k=1}^{N-1}{I_k \otimes I_k \over
(I_k,I_k)}$ for any orthogonal basis of the Cartan subalgebra.

 Consider the covariant algebra ${\Bbb C}[SL_N]_q$ as in section 5. It is
well known that
$${\Bbb C}[SL_N]_q \simeq \bigoplus_ \lambda{\rm End}(V_
\lambda)^*,\eqno(A1.8)$$
with $V_ \lambda$ being the simple $U_q \frak{sl}_N$-modules from the class
described in section 2. Hence the operator $R_{{\Bbb C}[SL_N]_q{\Bbb
C}[SL_N]_q}$ is well defined.

 We follow V. Drinfeld's approach in defining such $\widetilde{w}_0 \in{\Bbb
C}[SL_N]_q^*$ that
$$\langle f_1 \cdot f_2,\widetilde{w}_0 \rangle=\langle R_{{\Bbb
C}[SL_N]_q{\Bbb C}[SL_N]_q}(f_1 \otimes f_2),\widetilde{w}_0 \otimes
\widetilde{w}_0 \rangle \eqno(A1.9)$$
for all $f_1,f_2 \in{\Bbb C}[SL_N]_q$.

 Consider the Hopf algebra ${\Bbb C}[SL_2]_q$. Remind the relation
$t_{12}t_{21}=t_{21}t_{12}$. It is well known that every element of this
algebra admits a unique decomposition as follows
$$f = \sum_{j=1}^\infty t_{22}^j \cdot
f_j(t_{12},t_{21})+f_0(t_{12},t_{21})+\sum_{j=1}^\infty
f_{-j}(t_{12},t_{21})\cdot t_{11}^j,$$
with $f_j$ being polynomials in two commuting indeterminates.

 Consider the element $\overline{s}\in{\Bbb C}[SL_2]_q^*$ given by
$$\overline{s}(f)=f_0(q,-1),\qquad f \in{\Bbb C}[SL_2]_q$$
(it is a q-analogue of the Weyl element $\pmatrix{0 & 1 \cr -1 & 0}$).

 Associate to each $j=1,\ldots,N-1$ the onto homomorphism $\varphi_j:{\Bbb
C}[SL_N]_q \to{\Bbb C}[SL_2]_q$,
$$\qquad\varphi_j(t_{ik})=\left \{\begin{array}{ccl}\delta_{ik} &,& i \notin
\{j,j+1 \}\quad{\rm or}\quad k \notin \{j,j+1 \}\\ t_{i-j+1,k-j+1} &,& {\rm
otherwise}\end{array}\right..$$

 Consider a reduced decomposition $w_0=s_{i_1}s_{i_2}\ldots s_{i_M}$,
$M=N(N-1)/2$, of the longest permutation $w_0 \in S_N$, together with the
element
$$\overline{w}_0=\overline{s}_{i_1}\overline{s}_{i_2}\ldots
\overline{s}_{i_M}\in{\Bbb C}[SL_N]_q^*,\eqno(A1.10)$$
with $\overline{s}_j=\overline{s}\circ \varphi_j$, $j=1,\ldots,N-1$.It is
well known that $\overline{w}_0$ is independent of the choice of reduced
decomposition. Now we are in a position to define $\widetilde{w}_0$ by
$$\widetilde{w}_0=\overline{w}_0^{-1}\cdot q^{-{1 \over
2}\sum_kI_k^2/(I_k,I_k)},\eqno(A1.11)$$
with $\{I_k \}_{k=1}^{N-1}$ being an orthogonal basis of the Cartan
subalgebra. (It follows from \cite{LS2} that this is well defined and (A1.9)
holds.)

 Note that in \cite{LS1, LS2} the 'quantum simple maps' have been used to
produce automorphisms $T_i$ of $U_q \frak{sl}_N$ involved into the
definition of $E_{\beta_j}$, $F_{\beta_j}$. Specifically, one has
$$T_i(\xi)=\overline{s}_i \cdot \xi \cdot \overline{s}_i^{-1},\qquad \xi
\in U_q \frak{sl}_N,\qquad i=1,\ldots,N-1.$$
Hence $T_i$, $i=1,\ldots,N-1$, are extendable by a continuity from the
weakly dense subalgebra $U_q \frak{sl}_N \subset{\Bbb C}[SL_N]_q^*$ onto the
entire ${\Bbb C}[SL_N]_q^*$.

\bigskip

\section*{Appendix 2. On some $*$-representation of \boldmath ${\rm
Pol}(\widetilde{X})_q$}

\subsection*{\S 1. The construction of a $*$-representation \boldmath
$\widetilde{\prod}$}

 In section 6 ${\Bbb C}[SL_N]_q$ was equipped with involutions $*$ and
$\star$. In view of (\ref{*Tij}), (\ref{tijst}) one has
$$t_{ij}^*=\lambda_1(i)\lambda_2(j)t_{ij}^\star,\qquad
i,j=1,\ldots,N,\eqno(A2.1.1)$$
with
$$\lambda_1(k)={\rm sign}(k-m-1/2), \qquad \lambda_2(k)={\rm
sign}(n-k+1/2).\eqno(A2.1.2)$$

 A representation $\pi$ of ${\Bbb C}[SL_N]_q$ in a pre-Hilbert space
determines a $*$-representation of ${\rm Pol}(\widetilde{X})_q$ if and only
if $\pi(t_{ij})^*=\lambda_1(i)\lambda_2(j)\pi(t_{ij}^\star)$ for all
$i,j=1,\ldots,N$.

 Our purpose is to produce such a $*$-representation $\widetilde{\prod}$ of
${\rm Pol}(\widetilde{X})_q$ that $\widetilde{\prod}(x)\ne 0$. The method we
apply is well known in quantum group theory \cite{CP}.

 Let $\Lambda'=(\lambda'(1),\lambda'(2),\ldots,\lambda'(N))$,
$\Lambda''=(\lambda''(1),\lambda''(2),\ldots,\lambda''(N))$ be two sequences
whose entries are $\pm 1$. Suppose we are given a representation $\pi$ of
${\Bbb C}[SL_N]_q$ in a pre-Hilbert space. $\pi$ is said to be of type
$(\Lambda',\Lambda'')$ if
$$\pi(t_{ij})^*=\lambda'(i)\lambda''(j)\pi(t_{ij}^\star).$$
(In the special case of the sequences $(\Lambda_1,\Lambda_2)$ determined by
(A2.1.2) one has the class of all $*$-representations of ${\rm
Pol}(\widetilde{X})_q$).

\medskip

\noindent{\bf Proposition  A2.1.1.} {\it Suppose that the representations
$\pi'$ and $\pi''$ are of types $(\Lambda',\Lambda'')$ and
$(\Lambda'',\Lambda''')$ respectively. Then their tensor product
$\pi'\otimes \pi''$ is of type $(\Lambda',\Lambda''')$.}

\smallskip

 {\bf Proof.} An application of the relation $(\lambda''(k))^2=1$ and the
fact that the comultiplication $\Delta:{\Bbb C}[SU_N]_q \to{\Bbb
C}[SU_N]_q^{\otimes 2}$ is a homomorphism of $*$-algebras yields
$$(\pi'\otimes \pi''(t_{ij}))^*=\sum_{k=1}^N \pi'(t_{ik})^*\otimes
\pi''(t_{kj})^*=\lambda'(i)\lambda'''(j)\sum_{k=1}^N(\lambda''(k))^2
\pi'(t_{ik}^\star)\otimes \pi''(t_{kj}^\star)=$$
$$=\lambda'(i)\lambda'''(j)\pi'\otimes \pi''(t_{ij}^\star)\eqno \Box$$

\medskip

 {\sc Example A2.1.2.} In the special case $m=n=1$ one has
$\Lambda_1=(-1,1)$, $\Lambda_2=(1,-1)$, $t_{11}^*=-t_{11}^\star$,
$t_{12}^*=t_{12}^\star$, $t_{21}^*=t_{21}^\star$, $t_{22}^*=-t_{22}^\star$.
Let $\{e_j \}_{j \in{\Bbb Z}_+}$ be such an orthogonal basis of a
pre-Hilbert space that $(e_0,e_0)=1$,
$(e_j,e_j)=(q^{-2}-1)(q^{-4}-1)\ldots(q^{-2j}-1)$ for $j \in{\Bbb N}$. The
following formulae determine a representation $\pi_+$ of type
$(\Lambda_1,\Lambda_2)$:
$$\begin{array}{cc}\pi_+(t_{12})e_j=q^{-j}e_j, &
\pi_+(t_{21})e_j=-q^{-(j+1)}e_j, \\ \pi_+(t_{11})e_j=e_{j+1}, &
\pi_+(t_{22})e_j=(1-q^{-2j})e_{j-1}\end{array}.\eqno(A2.1.3)$$

\medskip

 {\sc Example A2.1.3.} Let $N \ge 2$, $k \in \{1,\ldots,N-1 \}$, and the
pair $(\Lambda',\Lambda'')$ possesses the properties:
$\lambda'(j)=\lambda''(j)$ for $j \notin \{k,k+1 \}$, $\lambda'(k)=-1$,
$\lambda''(k)=1$, $\lambda'(k+1)=1$, $\lambda''(k+1)=-1$. Consider the
homomorphism of algebras
$$\psi_k=\psi_{(k,k+1)}:{\Bbb C}[SL_N]_q \to{\Bbb C}[SL_2]_q,\qquad
\psi_k(t_{ij})=\left \{\begin{array}{ccl}t_{i-k+1,j-k+1} &,& i,j \in \{k,k+1
\}\\ \delta_{ij} &,&{\rm otherwise}\end{array}\right..$$ It is well known
that $\psi_k(f^\star)=(\psi_k(f))^\star$ for all $f \in{\Bbb C}[SL_N]_q$. On
can readily deduce from the definitions that the representation $\pi_+\circ
\psi_k$ of ${\Bbb C}[SL_N]_q$ is of type $(\Lambda',\Lambda'')$.

\medskip

 Now turn to a construction of a $*$-representation $\widetilde{\prod}$ of
${\rm Pol}(\widetilde{X})_q$, that is, a representation of ${\Bbb
C}[SL_N]_q$ of type $(\Lambda_1,\Lambda_2)$.

 Consider the element
$$\pmatrix{1 & 2 & \ldots & m & m+1 & m+2 & \ldots & N \cr n+1 & n+2& \ldots
& N & 1 & 2 & \ldots & n}$$
of the symmetric group $S_N$, together with its reduced decomposition
$\sigma_1 \cdot \sigma_2 \cdot \ldots \cdot \sigma_{mn}$. Let $s_0=e$,
$s_1=\sigma_1$, $s_2=\sigma_1 \cdot \sigma_2$, $\ldots$, $s_{mn}=\sigma_1
\cdot \sigma_2 \cdot \ldots \cdot \sigma_{mn}$. Our construction involves
the sequence $\Lambda^{(0)}$, $\Lambda^{(1)}$, $\ldots$, $\Lambda^{(mn)}$,
given by
$$\Lambda^{(j)}=\left(\lambda_2(s_{mn-j}(1)),\lambda_2(s_{mn-j}(2)),\ldots,
\lambda_2(s_{mn-j}(N))\right).$$

 Evidently, $\Lambda^{(0)}=\Lambda_1$, $\Lambda^{(mn)}=\Lambda_2$, and the
sequences in each pair $(\Lambda^{(j)},\Lambda^{(j+1)})$, $j=1,\ldots,mn-1$,
differ only by a permutation of some two neighbour terms $+1$, $-1$. Just as
in Example A2.1.3, construct representations of types
$(\Lambda^{(j)},\Lambda^{(j+1)})$. Their tensor product is of type
$(\Lambda_1,\Lambda_2)$ due to proposition A2.1.1. Hence this tensor product
$\widetilde{\prod}$ is a $*$-representation of ${\rm Pol}(\widetilde{X})_q$.

 We are interested in considering the restriction of $\widetilde{\prod}$
onto the subalgebra ${\rm Pol}(X)_q$ (see section 7).

\bigskip

\subsection*{\S 2. A faithful irreducible $*$-representation of \boldmath
${\rm Pol}(\widetilde{X})_q$}

\noindent{\bf Lemma A2.2.1.} {\it Let $v$ be such a vector in the space of a
$*$-representation $\rho$ of ${\rm Pol}(\widetilde{X})_q$ that
$${\rho(t_{\{m+1,m+2,\ldots,N \}\{1,2,\ldots,n \}}^{\wedge n})v=c \cdot
v,\qquad c \in{\Bbb C};\atop \rho(t_{\{m+1,m+2,\ldots,N \}J}^{\wedge
n})v=0,\qquad J \ne \{1,2,\ldots,n \}}.\eqno(A2.2.1)$$
Then $|c|=q^{mn}$.}

\smallskip

 {\bf Proof.} There is a relation in ${\rm Pol}(\widetilde{X})_q$ between
$t_{\{1,2,\ldots,m \}I}^{\wedge m}$ and $t_{\{m+1,m+2,\ldots,N \}J}^{\wedge
n}$ derived from ${\rm det}_qT=1$, $T=(t_{ij})_{i,j=1,\ldots,N}$, via a
q-analogue of the Laplace formula. By a virtue of (A2.2.1),
$$\rho \left((-q)^{mn}t_{\{1,2,\ldots,m \}\{n+1,n+2,\ldots,N \}}^{\wedge m}
\cdot t_{\{m+1,m+2,\ldots,N \}\{1,2,\ldots,n \}}^{\wedge n}\right)v=v.$$
On the other hand,
$$t_{\{1,2,\ldots,m \}\{n+1,n+2,\ldots,N \}}^{\wedge
m}=(-q)^{mn}\left(t_{\{m+1,m+2,\ldots,N \}\{1,2,\ldots,n \}}^{\wedge
n}\right)^*.$$
So
$$q^{2mn}\left \|\rho \left(t_{\{m+1,m+2,\ldots,N \}\{1,2,\ldots,n
\}}^{\wedge n}\right)v \right \|^2=\|v \|^2,$$
that is,
$$\left \|\rho \left(t_{\{m+1,m+2,\ldots,N \}\{1,2,\ldots,n \}}^{\wedge n}
\right)v \right \|=q^{-mn}\|v \|.\eqno \Box$$

\medskip

 Impose the notation
$${\bf e}_\Bbbk=e_{k_1}\otimes e_{k_2}\otimes \ldots \otimes
e_{k_{mn}},\qquad \Bbbk=(k_1,k_2,\ldots,k_{mn})\in{\Bbb Z}_+^{mn}$$
for basis vectors of the space of $\widetilde{\prod}$.

\medskip

 {\sc Example A2.2.2.} Let $m=n=2$. It follows from the definitions that
$$\psi_2 \otimes \psi_3 \otimes \psi_1 \otimes
\psi_2(t_{13}t_{24}-qt_{14}t_{23})=t_{12}\otimes t_{12}\otimes t_{12}\otimes
t_{12}.$$
(In fact,
$$\psi_2 \otimes \psi_3 \otimes \psi_1 \otimes \psi_2(t_{13})=1 \otimes 1
\otimes t_{12}\otimes t_{12},$$
$$\psi_2 \otimes \psi_3 \otimes \psi_1 \otimes \psi_2(t_{24})=t_{12}\otimes
t_{12}\otimes 1 \otimes 1,$$
$$\psi_2 \otimes \psi_3 \otimes \psi_1 \otimes \psi_2(t_{14})=0.$$
Hence for all $\Bbbk=(k_1,k_2,k_3,k_4)$
$$\widetilde{\prod}(t_{13}t_{24}-qt_{14}t_{23}){\bf
e}_\Bbbk=q^{-(k_1+k_2+k_3+k_4)}{\bf e}_\Bbbk.\eqno(A2.2.2)$$

\medskip

 It is easy to extend (A2.2.2) onto the case of arbitrary $m,n \in{\Bbb N}$.

 Consider the element $u=(m+1,m+2,\ldots,N,1,2,\ldots,m)\in S_N$, together
with its reduced decomposition of the form $u=\sigma_1 \sigma_2 \sigma_3
\ldots \sigma_{mn}$,$$\sigma_k=\left(m-\left[{\textstyle k-1 \over
\textstyle n}\right]+\left \{{\textstyle k-1 \over \textstyle n}\right
\}n,m-\left[{\textstyle k-1 \over \textstyle n}\right]+\left
\{{\textstyle k-1 \over \textstyle n}\right \}n+1 \right)$$ (here $[\cdot]$,
$\{\cdot \}$ stand for integral and fractional parts of a real number,
respectively). For example, in the case $m=2$, $n=3$, one has
$u=(3,4,5,1,2)$, and the above reduced decomposition acquires the form
$u=(2,3)(3,4)(4,5)(1,2)(2,3)(3,4)$.

 It is easy to show that $\widetilde{\prod}\pi_+^{\otimes mn}\circ \Psi$,
with $\Psi:{\Bbb C}[SL_N]_q \to{\Bbb C}[SL_2]_q^{\otimes mn}$,
$\Psi=\psi_{\sigma_1}\otimes \psi_{\sigma_2}\otimes \ldots \otimes
\psi_{\sigma_{mn}}$ (we use here the notation from Example(A2.1.3)).

\medskip

\noindent{\bf Lemma A2.2.2.} {\it For all $\Bbbk \in{\Bbb Z}_+^{mn}$
$$\widetilde{\prod}\left(t_{\{1,\ldots,m
\}\{n+1,\ldots,N \}}\right){\bf e}_\Bbbk=q^{-\sum \limits_jk_j}{\bf
e}_\Bbbk.\eqno(A2.2.3)$$}

\smallskip

 {\bf Proof.} Let $n$ be fixed. We use an induction in $m$ to show that for
$i \le m$
$$\Psi(=\Psi_m):t_{ij}\mapsto \left \{\begin{array}{ccl}0 &,& j>i+n \\
\underbrace{1 \otimes \ldots \otimes 1}_{(m-i)n}\otimes
\underbrace{t_{12}\otimes \ldots \otimes t_{12}}_n \otimes \underbrace{1
\otimes \ldots \otimes 1}_{(i-1)n}&,& j=i+n \end{array}\right..$$

 In the case $m=1$ the statement is evident (since $\Psi_1(t_{1,n+1})=\psi_1
\otimes \psi_2 \otimes \ldots \otimes \psi_n \left(\sum t_{1i_1}\otimes
t_{i_1i_2}\otimes \ldots \otimes t_{i_{n-1},n+1}\right)=\psi_1 \otimes \psi_2
\otimes \ldots \otimes \psi_n(t_{12}\otimes t_{23}\otimes \ldots \otimes
t_{n,n+1})=t_{12}\otimes t_{12}\otimes \ldots \otimes t_{12})$.

 Now we are to make the induction passage from $(m-1)$ to $m$. Let
$\Phi_m=\psi_{\sigma_1}\otimes \psi_{\sigma_2}\otimes \ldots \otimes
\psi_{\sigma_n}:{\Bbb C}[SL_N]_q \to{\Bbb C}[SL_2]_q^{\otimes n}$,
$\Psi_{m-1}'=\psi_{\sigma_{n+1}}\otimes \psi_{\sigma_{n+2}}\otimes \ldots
\otimes \psi_{\sigma_{mn}}:{\Bbb C}[SL_N]_q \to{\Bbb C}[SL_2]_q^{\otimes
n(m-1)}$, i.e. $\Psi_m=\Phi_m \otimes \Psi_{m-1}'$.

 Obviously, the subalgebra of ${\Bbb C}[SL_N]_q$ generated by $t_{ij}$ with
$i,j<N$, is isomorphic to ${\Bbb C}[SL_{N-1}]_q$. If we agree to identify in
what follows ${\Bbb C}[SL_{N-1}]_q$ with this subalgebra, one can claim that
$$\Psi_{m-1}(t_{ij})=\Psi_{m-1}'(t_{ij}),\qquad (i,j<N).$$

 Besides that,
$\Psi_{m-1}'(t_{iN})=\Psi_{m-1}'(t_{Ni})=\delta_{iN}\underbrace{1 \otimes
1 \otimes \ldots \otimes 1}_{(m-1)n}$. These facts are to be used in the
passage from $(m-1)$ to $m$. We start this passage with considering the
special case of elements of the last column:
$$\Psi_m(t_{iN})=\Phi_m \otimes \Psi_{m-1}'\left(\sum_{k=1}^Nt_{ik}\otimes
t_{kN}\right)=\sum_{k=1}^N \Phi_m(t_{ik})\otimes \Psi_{m-1}'(t_{kN})=$$
$$=\Phi_m(t_{iN})\otimes \Psi_{m-1}'(t_{NN})=\Phi_m(t_{iN})\otimes
\underbrace{1 \otimes \ldots \otimes 1}_{(m-1)n}.$$

 In the case $i<m$ it is easily deducible from the definition of $\Phi_m$
that $\Phi_m(t_{iN})=0$ (the cycle $\sigma_1 \sigma_2 \ldots
\sigma_n=(m,m+1)(m+1,m+2)\ldots(N-1,N)$ can not send $N$ to $i$). If $i=m$,
then
$$\Phi_m(t_{mN})=\psi_{\sigma_1}\otimes \psi_{\sigma_2}\otimes \ldots
\otimes \psi_{\sigma_n}(t_{mN})=$$
$$=\psi_m \otimes \psi_{m+1}\otimes \ldots \otimes \psi_{N-1}\left(\sum
t_{mi_1}\otimes t_{i_1i_2}\otimes \ldots \otimes t_{i_{n-1}N}\right)=$$
$$=\psi_m \otimes \psi_{m+1}\otimes \ldots \otimes
\psi_{N-1}(t_{m,m+1}\otimes t_{m+1,m+2}\otimes \ldots \otimes
t_{N-1,N})=\underbrace{t_{12}\otimes t_{12}\otimes \ldots \otimes
t_{12}}_n.$$

 Thus we have done an induction passage for elements of the last column.
What remains is to consider the case $j \le N-1$, $i+n \le j$ (note that
$i+n \le N-1$ implies $i<m$). One has
$$\Psi_m(t_{ij})=\Phi_m \otimes \Psi_{m-1}'\left(\sum_{k=1}^Nt_{ik}\otimes
t_{kj}\right),$$
and, since $\Psi_{m-1}'(t_{Nj})=0$ for $j<N$,
$$\Psi_m(t_{ij})=\sum_{k=1}^{N-1}\Phi_m(t_{ik})\otimes
\Psi_{m-1}'(t_{kj}).\eqno(*)$$

 Consider the element $\Phi_m(t_{ik})$:
$$\Phi_m(t_{ik})=\psi_m \otimes \psi_{m+1}\otimes \ldots \otimes
\psi_{N-1}\left(\sum t_{ij_1}\otimes t_{j_1j_2}\otimes \ldots \otimes
t_{j_{n-1}k}\right).$$

 Since $i<m$, one has $\psi_m(t_{ij_1})\ne 0$ only for $j_1=i$. Similarly,
$j_2=j_3=\ldots=k=i$. Hence, in (*) only one term "survives":
$$\Psi_m(t_{ij})=\Phi_m(t_{ii})\otimes \Psi_{m-1}(t_{ij}).$$

 The induction hypothesis implies $\Psi_{m-1}'(t_{ij})=0$ for $j>i+n$, and
hence for such $i$ and $j$ that $\Psi_m(t_{ij})=0$. What remains is to
consider the case $i+n=j$. If so, again the induction hypothesis yields
$$\Psi_{m-1}(t_{ij})=\underbrace{1 \otimes 1 \otimes \ldots \otimes
1}_{(m-1-i)n}\otimes \underbrace{t_{12}\otimes \ldots \otimes t_{12}}_n
\otimes \underbrace{1 \otimes 1 \otimes \ldots \otimes 1}_{(i-1)n},$$
i. e.
$$\Phi_m(t_{ii})\otimes \Psi_{m-1}'(t_{ij})=$$
$$=\underbrace{1 \otimes 1 \otimes \ldots \otimes 1}_n \otimes \underbrace{1
\otimes 1 \otimes \ldots \otimes 1}_{(m-1-i)n}\otimes
\underbrace{t_{12}\otimes \ldots \otimes t_{12}}_n \otimes \underbrace{1
\otimes 1 \otimes \ldots \otimes 1}_{(i-1)n}=$$
$$=\underbrace{1 \otimes 1 \otimes \ldots \otimes 1}_{(m-i)n}\otimes
\underbrace{t_{12}\otimes \ldots \otimes t_{12}}_n \otimes \underbrace{1
\otimes 1 \otimes \ldots \otimes 1}_{(i-1)n}.$$

 This completes the induction passage. \hfill $\Box$

\medskip

 By a virtue of (A2.2.3) the operator $\widetilde{\Pi}(x)$ is invertible,
and hence the representation $\widetilde{\Pi}$ admits a unique extension
onto the $*$-algebra ${\rm Pol}(\widetilde{X})_{q,x}$. Let
$\Pi=\widetilde{\Pi}\circ i$ be the $*$-representation of ${\rm Pol}({\rm
Mat}_{mn})_q$ deduced from the $*$ -homomorphism ${\cal I}:{\rm Pol}({\rm
Mat}_{mn})_q \to{\rm Pol}(\widetilde{X})_{q,x}$ described in section 6.
Equip the pre-Hilbert representation space $\widetilde{H}$ of
$\widetilde{\Pi}$ and the algebras ${\rm Pol}(\widetilde{X})_q$, ${\rm
Pol}({\rm Mat}_{mn})_q$ with the gradations:
$$\widetilde{H}=\bigoplus_{j=0}^\infty \widetilde{H}_j,\qquad
\widetilde{H}_j=\{v \in \widetilde{H}|\;\widetilde{\Pi}(x)v=q^{-2j}v \},$$
$${\rm deg}(t_{ij})=\left \{\begin{array}{ccl}1 &,& i \le m \;\&\;j \le n \\
-1 &,& i>m \;\&\;j>n \\ 0 &,& otherwise \end{array}\right..$$ $${\rm
deg}(z_a^\alpha)=1,\;{\rm deg}(z_a^\alpha)^*=-1$$

 It follows from the commutation relations (\ref{tijt_}), (\ref{tijt*}) that
$\widetilde{\Pi}$ allows one to equip $\widetilde{H}$ with the structure of
a graded ${\rm Pol}(\widetilde{X})_q$-module, and the representation $\Pi$
with a structure of a graded ${\rm Pol}({\rm Mat}_{mn})_q$-module.

\medskip

\noindent{\bf Proposition A2.2.3.} {\it The {\sl graded}
${\rm Pol}(X)_q$-module $\widetilde{H}$ is simple.\footnote{That is, it has
no nontrivial graded submodules}}

\smallskip

 {\bf Proof.} Let $L \subset \widetilde{H}$ be a nontrivial graded
submodule. It follows from lemma A2.2.2 that the operators
$\widetilde{\Pi}(t)$ and $\widetilde{\Pi}(t^*)$ are the same. Also,
$\widetilde{\Pi}(x)L \subset L$ implies $\widetilde{\Pi}(t)L \subset L$,
$\widetilde{\Pi}(t^*)L \subset L$. Hence $\widetilde{\Pi}(f)L \subset L$ for
all $f \in{\rm Pol}(\widehat{X})_q$ by lemma 8.2. In particular,
$t_{\{1,2,\ldots,m \}J}^{\wedge m}L \subset L$ for all $m$-element subsets
$J \subset \{1,2,\ldots,N \}$. On the other hand, in this case there exists
a non-zero vector $v \in L$ such that $\widetilde{\Pi}(t_{\{m+1,\ldots,N
\}I})v=0$ for all $n$-element subsets $I \subset \{1,2,\ldots,N \}$
different from $\{1,2,\ldots,n \}$. By a virtue of lemmas A2.2.1 and A2.2.2,
the subspace of all such vectors is one-dimensional. Hence, $v={\tt
const}\cdot{\bf e}_0$, ${\bf e}_0 \in L$, $L=\widetilde{H}$. A
contradiction. \hfill $\Box$

\medskip

 We prove in section 8 that the restriction of $\widetilde{\Pi}$ onto ${\rm
Pol}(\widehat{X})_q$ is a faithful representation of this subalgebra.

\twocolumn

\begin{theindex}

  \item $A$-module algebra, 2
  \item $A^{\rm op}$-module coalgebra, 2
  \item $J_{a \alpha}$, 16
  \item $U_q \frak{su}_N$, 21
  \item $U_q \frak{su}_n \otimes U_q \frak{su}_m$, 11
  \item $U_q \frak{su}_{nm}$, 10
  \item ${\Bbb C}[SU_N]_q$, 21
  \item ${\rm Pol}(\widetilde{X})_q$, 19

  \indexspace

  \item antimodule, 11

  \indexspace

  \item classical universal enveloping algebra, 2
  \item covariant algebra, 2
  \item covariant coalgebra, 2

  \indexspace

  \item differential algebra $\wedge({\rm Mat}_{mn})_q$, 10

  \indexspace

  \item finite functions, 24
  \item full differential calculus, 9

  \indexspace

  \item generalized Verma module, 3
  \item grading for $U_q \frak{sl}_N$-modules, 2

  \indexspace

  \item invariant integral, 32

  \indexspace

  \item localization, 15

  \indexspace

  \item order one differential calculus, 9

  \indexspace

  \item polynomial algebra on the quantum principal homogeneous space,
          19

  \indexspace

  \item q-determinant, 14
  \item quantum universal enveloping algebra, 1

\end{theindex}

\end{document}